\RequirePackage{fix-cm}
\documentclass[twocolumn]{svjour3}          
\smartqed  
\usepackage{amsmath,amssymb}
\usepackage{graphicx}
\usepackage{mathtools}
\usepackage{tikz}
\usepackage{subfigure}
\usepackage{mathptmx}      
\usepackage{latexsym}
\newtheorem{assumption}[theorem]{Assumption}
\begin{document}

\title{On Darcy-and Brinkman-Type Models for Two-Phase Flow in Asymptotically Flat Domains}



\author{Alaa Armiti-Juber \and Christian Rohde }


\institute{Alaa Armiti-Juber/ Christian Rohde \\
              \email{a.armiti/ crohde@mathematik.uni-stuttgart.de}      \\     
\at
              Institute for Applied Analysis and Numerical Simulation, University of Stuttgart, Pfaffenwaldring 57, 70569 Stuttgart, Germany 
}
\date{}

\maketitle

\begin{abstract}
We study two-phase flow for Darcy and Brinkman regimes. To reduce the computational complexity for flow in vertical equilibrium various simplified models have been suggested. Examples are dimensional reduction by vertical integration, the multiscale model approach in [Guo et al., 2014] or the asymptotic approach in [Yortsos, 1995]. The latter approach uses a geometrical scaling. We show the efficiency of the approach in asymptotically flat domains for Darcy regimes. Moreover, we prove that it is vastly equivalent to the multiscale model approach.

We apply then asymptotic analysis to the two-phase flow model in Brinkman regimes. The limit model is a single nonlocal evolution law with a pseudo-parabolic extension. Its computational efficiency is demonstrated using numerical examples. Finally, we show that the new limit model exhibits overshoot behaviour as it has been observed for dynamical capillarity laws [Hassanizadeh and Gray, 1993].


\end{abstract}

\section{Introduction}
\label{intro}
In this  paper we consider incompressible two-phase flow through a porous medium. If capillary forces are neglected two-phase flow processes are governed by hyperbolic equations for saturation and elliptic Darcy-like laws for the pressure \cite{Bear1972,Helmig1997}, resulting in 
a strongly coupled nonlinear system of differential equations. Standard discretization techniques (based on e.g.~the implicit pressure and explicit saturation (IMPES) treatment) are computationally quite expensive in spatially three-dimensional settings \cite{HuberHelmig1999}. Therefore it is an important issue to perform model reductions if possible, in particular for any kinds of limit regimes.
\medskip

Here, we focus on porous media in asymptotically flat domains with almost horizontal fluid flow. If the fluids in the domain segregate fast, their pressures can be assumed to be essentially hydrostatic. This assumption is called the vertical equilibrium assumption (the name of this assumption should not be mixed up with the Vertical Equilibrium model presented later). It usually applies to asymptotically flat domains as well as to domains with very high vertical permeability where vertical forces, like gravity, equilibrate very fast. 
\medskip

The vertical equilibrium assumption has been utilized by many researchers to integrate the three-dimensional two-phase flow model over the vertical coordinate and consequently to reduce the computational complexity of the model by dimensional reduction. For example, we refer to \cite{Lake1989} in the field of petroleum studies and to \cite{Bear1972,Helmig1997} in the field of hydrogeology, where the assumption is well-known under the name of Dupuit's assumption. Recently, this approach has received attention in the field of $\text{CO}_2$-sequestration in saline aquifers \cite{Becker,Gasda2009,Gasda2011}.
\medskip

Such vertically integrated models have the tendency to overestimate horizontal spreading speeds, mainly in domains where the vertical equilibrium assumption is not fully satisfied \cite{court2012}. Examples are heterogeneous porous media with low vertical permeability, media with relatively large thickness and cases when the flowing fluids have a small difference in their densities. Therefore, different approaches that modify the vertical integration approach have been developed. For example, a Multiscale model that relaxes the vertical equilibrium assumption is proposed in \cite{Guo2014}. This model consists of a coarse scale vertically integrated equation for the vertically integrated pressure in the horizontal plane and a set of fine scale one-dimensional equations for saturation in the vertical direction. The two equations are coupled via an operator that reconstructs the coarse scale pressure along the vertical coordinate. Another approach for multi-layered flat domains is proposed in \cite{Qin2014}. In this approach, the dynamics in each layer is described by a macroscale thickness-integrated mass balance equation coupled by appropriate exchange terms.
\medskip

In this paper we exploit a model reduction approach of Yortsos \cite{Yortsos} that relies on asymptotic analysis in terms of the 
height$-$length ratio of the domain. The reduced two-phase model is still posed in the fully-dimensional domain but consists of a single nonlocal equation for the saturation alone. We consider this approach for two different two-phase flow scenarios:
\begin{itemize}
 \item Asymptotically flat domains with the total velocity satisfying a Darcy-type law.
 \item Asymptotically flat domains with the total velocity obeying Brinkman's equations.
\end{itemize}
\medskip

As mentioned above this approach has been explored by Yortsos for a two-phase Darcy-type flow model \cite{Yortsos}, motivated by the numerical results in \cite{Coats1971,YokoyamaLake1981,ZaptaLake1981}. He called this model the Vertical Equilibrium model. This model has been further investigated to develop selection principles for finding upper bounds on the speed of the mixing zone in the case of miscible displacement \cite{MenonOtto2005,YortsosSalin2006}. Existence of weak solutions for a regularization of the model in the framework of functions of bounded variations is proved in \cite{Armiti-Juber2014}. For miscible displacement in the gap of a planar Hele-Shaw or a long cylindrical capillary, the validity of this model is investigated in \cite{Yang1997}. A similar model is developed for layered reservoirs in \cite{Zhang2011,Zhang2012}. For media with parallel layers of contrasting porosity, absolute permeability and relative permeability, a generalization of the model is proposed in \cite{Debbabi2017}.
\medskip

We first suggest a conservative and stable finite-volume scheme for the nonlocal Vertical Equilibrium model. We show by numerical experiments that the Vertical Equilibrium model is much more accurate than vertically integrated models while still being computationally much more efficient
than the full two-phase flow model. In addition to this, we show that the Vertical Equilibrium model and the Multiscale model proposed in \cite{Guo2014} coincide in general. This surprising result might lead to further reduction in the computational complexity of the Multiscale model.  
\medskip 

In flat domains, we apply asymptotic analysis to the two-phase flow model with the second-order Brinkman's equations following \cite{Yortsos}. By doing this, we derive a nonlocal evolution equation for saturation only. It involves additionally a third-order term with  mixed --time and spatial-- derivatives. We propose a new numerical scheme and show again by numerical examples that this Brinkman Vertical Equilibrium model is in excellent agreement with the full two-phase flow model in flat domains, while being computationally much more efficient.
\medskip

The Brinkman Vertical Equilibrium model resembles the third-order pseudo-parabolic model proposed in \cite{HG1993} and mathematically analyzed in \cite{FanPop2010,vanDuijn2013}. This model has been derived in a completely other context to model rate-dependent capillarity pressures. In particular, it has been shown (see e.g.~\cite{vanDuijn2013}) that it supports the upset of overshooting of invading wetting fronts. We also observe overshooting in the Brinkman Vertical Equilibrium model \emph{without} introducing additional second gradient forces like in \cite{HG1993}.
\medskip

We conclude the introduction with an overview of the paper's content: In Section \ref{sec:background} we review the models that describe two-phase flow in asymptotically flat domains. Section \ref{sec:VE-numericalresults} provides several comparison examples of the Vertical Equilibrium model with the other models presented in the previous section. In Section \ref{sec:Brinkmanmodel}, we derive the Brinkman Vertical Equilibrium model to describe two-phase flow in flat domains. Section \ref{sec:B-numericalresults} provides a set of numerical tests that show the accuracy and the computational efficiency of the Brinkman Vertical Equilibrium model and its ability to support overshooting fronts.

\section{Models for Two-Phase Darcy Flow in Asymptotically Flat Domains}
\label{sec:background}
To study fluid flow in porous media on the macroscale physical properties of the phases are averaged over an appropriate control volume, the representative elementary volume \cite{Helmig1997}. On this scale, the averaged velocity of the fluids can be described by Darcy's law. 

In this section, we present first the governing equations for a two-phase flow model on the macroscale in a bounded domain $\Omega\subset\mathbb{R}^3$. Then, we review three different models for asymptotically flat domains: a Vertically Integrated model, the Multiscale model proposed by Guo \emph{et al.} \cite{Guo2014} and the Vertical Equilibrium model proposed by Yortsos \cite{Yortsos}.

\begin{figure}
\centering
\begin{tikzpicture}[scale = 0.75,>=latex]
\draw[dashed, thick] (0,0) -- coordinate (y axis mid) (0,3);
\node[rotate=90, above=0.3cm, black!40!blue] at (y axis mid) {$\partial\Omega_{\text{inflow}}$};
\node[rotate=90, above=-7.4cm, black!60!orange] at (y axis mid) {$\partial\Omega_{\text{outflow}}$};
\draw[dashed, thick] (9,0) -- coordinate (x axis mid) (9,3);
\node[left=4cm , above=1.3cm] at (x axis mid) {$\partial\Omega_{\text{imp}}$};
\node[left=4cm , above=-1.9cm] at (x axis mid) {$\partial\Omega_{\text{imp}}$};
\draw[very thick] (0,0) -- (9,0);
\draw[very thick] (0,3) -- (9,3); 

\draw[|-] (0,-1)--(4,-1) ;\draw[-|] (5,-1)--(9,-1);
\draw[very thick] (4.3,-1) node[anchor=west] {L};
\draw[|-] (-1,0)--(-1,1.3) ;\draw[-|] (-1,1.8)--(-1,3);
\draw[very thick] (-1,1.25) node[anchor=south] {H};

\draw[->,blue](-0.3,0.2)-- (-0.1,0.2);\draw[->,blue](-0.3,0.4)-- (-0.1,0.4);\draw[->,blue](-0.3,0.6)-- (-0.1,0.6);\draw[->,blue](-0.3,0.8)-- (-0.1,0.8);\draw[->,blue](-0.3,1)-- (-0.1,1);\draw[->,blue](-0.3,1.2)-- (-0.1,1.2);\draw[->,blue](-0.3,1.4)-- (-0.1,1.4);\draw[->,blue](-0.3,1.6)-- (-0.1,1.6);\draw[->,blue](-0.3,1.8)-- (-0.1,1.8);\draw[->,blue](-0.3,2)-- (-0.1,2);\draw[->,blue](-0.3,2.2)-- (-0.1,2.2);\draw[->,blue](-0.3,2.4)-- (-0.1,2.4);\draw[->,blue](-0.3,2.6)-- (-0.1,2.6);\draw[->,blue](-0.3,2.8)-- (-0.1,2.8);

\draw[->,black!60!orange](9.1,0.2)-- (9.3,0.2);\draw[->,black!60!orange](9.1,0.4)-- (9.3,0.4);\draw[->,black!60!orange](9.1,0.6)-- (9.3,0.6);\draw[->,black!60!orange](9.1,0.8)-- (9.3,0.8);\draw[->,black!60!orange](9.1,1)-- (9.3,1);\draw[->,black!60!orange](9.1,1.2)-- (9.3,1.2);\draw[->,black!60!orange](9.1,1.4)-- (9.3,1.4);\draw[->,black!60!orange](9.1,1.6)-- (9.3,1.6);\draw[->,black!60!orange](9.1,1.8)-- (9.3,1.8);\draw[->,black!60!orange](9.1,2)-- (9.3,2);\draw[->,black!60!orange](9.1,2.2)-- (9.3,2.2);\draw[->,black!60!black!60!orange](9.1,2.4)-- (9.3,2.4);\draw[->,black!60!orange](9.1,2.6)-- (9.3,2.6);\draw[->,black!60!orange](9.1,2.8)-- (9.3,2.8);

\draw[-] (0,0) -- (0.1,-0.2);\draw[-] (0.3,0) -- (0.4,-0.2);\draw[-] (0.6,0) -- (0.7,-0.2);\draw[-] (0.9,0) -- (1,-0.2);\draw[-] (1.2,0) -- (1.3,-0.2);\draw[-] (1.5,0) -- (1.6,-0.2);\draw[-] (1.8,0) -- (1.9,-0.2);\draw[-] (2.1,0) -- (2.2,-0.2);\draw[-] (2.4,0) -- (2.5,-0.2);\draw[-] (2.7,0) -- (2.8,-0.2);\draw[-] (3,0) -- (3.1,-0.2);\draw[-] (3.3,0) -- (3.4,-0.2);\draw[-] (3.6,0) -- (3.7,-0.2);\draw[-] (3.9,0) -- (4,-0.2);\draw[-] (4.2,0) -- (4.3,-0.2);\draw[-] (4.5,0) -- (4.6,-0.2);\draw[-] (4.8,0) -- (4.9,-0.2);\draw[-] (5.1,0) -- (5.2,-0.2);\draw[-] (5.4,0) -- (5.5,-0.2);\draw[-] (5.7,0) -- (5.8,-0.2);\draw[-] (6,0) -- (6.1,-0.2);\draw[-] (6.3,0) -- (6.4,-0.2);\draw[-] (6.6,0) -- (6.7,-0.2);\draw[-] (6.9,0) -- (7,-0.2);\draw[-] (7.2,0) -- (7.3,-0.2);\draw[-] (7.5,0) -- (7.6,-0.2);\draw[-] (7.8,0) -- (7.9,-0.2);\draw[-] (8.1,0) -- (8.2,-0.2);\draw[-] (8.4,0) -- (8.5,-0.2);\draw[-] (8.7,0) -- (8.8,-0.2);\draw[-] (9,0) -- (9.1,-0.2);

\draw[-] (0,3) -- (0.1,3.2);\draw[-] (0.3,3) -- (0.4,3.2);\draw[-] (0.6,3) -- (0.7,3.2);\draw[-] (0.9,3) -- (1,3.2);\draw[-] (1.2,3) -- (1.3,3.2);\draw[-] (1.5,3) -- (1.6,3.2);\draw[-] (1.8,3) -- (1.9,3.2);\draw[-] (2.1,3) -- (2.2,3.2);\draw[-] (2.4,3) -- (2.5,3.2);\draw[-] (2.7,3) -- (2.8,3.2);\draw[-] (3,3) -- (3.1,3.2);\draw[-] (3.3,3) -- (3.4,3.2);\draw[-] (3.6,3) -- (3.7,3.2);\draw[-] (3.9,3) -- (4,3.2);\draw[-] (4.2,3) -- (4.3,3.2);\draw[-] (4.5,3) -- (4.6,3.2);\draw[-] (4.8,3) -- (4.9,3.2);\draw[-] (5.1,3) -- (5.2,3.2);\draw[-] (5.4,3) -- (5.5,3.2);\draw[-] (5.7,3) -- (5.8,3.2);\draw[-] (6,3) -- (6.1,3.2);\draw[-] (6.3,3) -- (6.4,3.2);\draw[-] (6.6,3) -- (6.7,3.2);\draw[-] (6.9,3) -- (7,3.2);\draw[-] (7.2,3) -- (7.3,3.2);\draw[-] (7.5,3) -- (7.6,3.2);\draw[-] (7.8,3) -- (7.9,3.2);\draw[-] (8.1,3) -- (8.2,3.2);\draw[-] (8.4,3) -- (8.5,3.2);\draw[-] (8.7,3) -- (8.8,3.2);\draw[-] (9,3) -- (9.1,3.2);

\draw[fill=blue!50] (0,3)-- (0,0) -- (1.5,0)[rounded corners=10pt]--(2.5,0.5)--(1.8,1)--(3,2.2)--(2,3);
\draw[fill=orange!30]  (9,3)--(9,0)-- (1.5,0)[rounded corners=10pt]--(2.5,0.5)--(1.8,1)--(3,2.2) -- (2,3);
\draw[thick] (8,2.5)  node[anchor=west] {$\Omega_v$}; 
\end{tikzpicture}
\caption{An illustration of the displacement process in a vertical-cross section $\Omega_v$ of a three-dimensional domain $\Omega$. For an asymptotically flat domain, the ratio $H/L$ tends to zero.}
\label{fig:displacementprocess}
\end{figure}
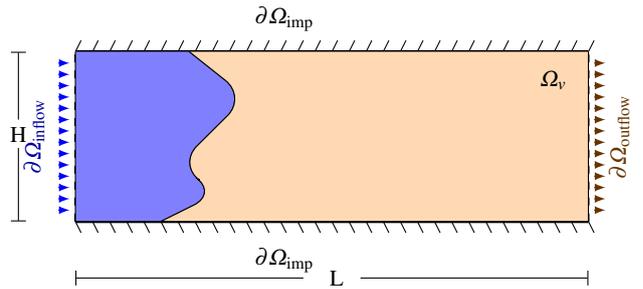
\subsection{The Two-Phase Flow Model in the Fully-Dimensional Domain $\Omega$}
\label{sec:two-phase2}
In the saturated case, the two-phase flow model (TP-model) is a combination of the continuity equation, Darcy's law, and the incompressibility relation \cite{Helmig1997}, i.e.,
\begin{equation}\left.
\begin{array}{rl}
\phi\partial_{t}S_{\alpha}+\nabla\cdot\textbf{v}_{\alpha}=&0,\\
\textbf{v}_{\alpha}=&-\lambda_{\alpha}(S_{\alpha})\textbf{K}\left(\nabla p_{\alpha}- \rho_{\alpha}g\textbf{e}_z\right),\\
\nabla\cdot\textbf{v}=&0
\end{array}\right.
\label{eq:two-phase}
\end{equation}
in $\Omega\times (0,T)$ for both, the invading $\alpha=i$ and the defending $\alpha=d$ phase. The unknowns in \eqref{eq:two-phase} are the saturation $S_{i}=S_{i}(x,y,z,t) \in [0,1]$ and the pressure $p_{i}=p_{i}(x,y,z,t)\in \mathbb{R}$. Saturation and pressure of the defending phase are then computed using the closure relations $S_i+S_d=1$ and $p_d-p_i=p_c$, respectively, where the capillary pressure $p_c$ can be assumed to be a given function of saturation of the invading phase \cite{Genuchten1980}, $p_c=p_c(S_i)$. The total velocity $\textbf{v}=\textbf{v}(x,y,z,t)\in\mathbb{R}^{3}$, for any $(x,y,z,t)\in \Omega\times[0,T]$, is the sum of the phases' velocities $\textbf{v}=\textbf{v}_i+\textbf{v}_d$. It consists of two horizontal components $(u,v)$ and one vertical component $(w)$, i.e., $\textbf{v}=(u,v,w)^T$. The intrinsic permeability tensor 
\begin{align}
 \textbf{K} = \textbf{K}(x,y,z) = \left(\begin{array}{c c c}
               K_x(x,y,z) & 0 & 0\\
               0 & K_y(x,y,z) & 0\\
               0 & 0 & K_z(x,y,z)
              \end{array}
\right),
\label{eq:permeabilitytensor}
\end{align}
is given with $K_x,\,K_y$ and $K_z$ being the permeabilities in $x$-, $y$-, and $z$-directions. $\phi=\phi(x,y,z)\in [0,1]$ is the medium's porosity, $\rho_{\alpha}$ is the phase density, $g$ is the gravitational acceleration, and $\textbf{e}_z=(0,0,1)^T$. The phase mobility $\lambda_{\alpha}:[0,1]\rightarrow [0,\infty)$ is given by $\lambda_{\alpha}(S_{\alpha}):=k_{r\alpha}(S_{\alpha})/\mu_{\alpha}$, where $k_{r\alpha}$ is the phase relative permeability and the constant $\mu_{\alpha}>0$ is the phase viscosity. 

For later purpose, we introduce the sum of the phases' mobilities, called the total mobility $\lambda_{tot}=\lambda_i + \lambda_d$. The fractional flow function for the invading phase is defined as $f = \tfrac{\lambda_i}{\lambda_{tot}}$. Using quadratic relative permeabilities results in the specific choice
 	\begin{align}
 	 f(S_i)= \dfrac{M S_i^2}{M S_i^2+ (1-S_i)^2}, \quad \lambda_{tot}(S_i)=M S_i^2+ (1-S_i)^2.
 	\end{align}
Here $M\coloneqq \mu_d /\mu_i$ is the viscosity ratio of the defending phase and the invading phase.

The domain $\Omega\subset\mathbb{R}^3$ is initially saturated with a resident saturation $S_d$. The boundary $\partial \Omega$, as illustrated in Figure \ref{fig:displacementprocess}, is decomposed into an inflow boundary $\partial\Omega_{\text{inflow}}$, an outflow boundary $\partial\Omega_{\text{inflow}}$ and an impermeable boundary $\partial \Omega_{\text{imp}}$. The initial and boundary conditions are summarized as follows
\begin{align}
\begin{array}{rll}
 S_i(\cdot,\cdot,\cdot,0)&=0, \quad\quad &\text{ in } \Omega, \\
 S_i&=S_{\text{inflow}},\quad\quad &\text{ on } \partial\Omega_{\text{inflow}}\times [0,T],\\
 \textbf{n}\cdot\textbf{v}_{\alpha}&= 0, \quad\quad &\text{ on }\partial \Omega_{\text{imp}}\times [0,T],
\end{array}
\label{eq:TPF-IBC}
\end{align}
where $S_{\text{inflow}}\in [0,1]$ and $\textbf{n}$ denotes the outer normal vector of the boundary $\partial \Omega_{\text{imp}}$. Throughout the paper, the following assumptions are imposed:
\begin{assumption}
\begin{enumerate}
\item The permeability components $K_{x}=K_y=K_z$ in \eqref{eq:permeabilitytensor} are strictly positive and bounded.
\item There exists a constant $c>0$ s.t. $\int_{0}^{H} \phi (\cdot,\cdot,z)\,dz>c$. 
\item For the total mobility $\lambda_{tot}= \lambda_{tot}(S_i) \in C^2([0,1])$ there is a constant $a>0$ with $0<a<\lambda_{tot}$.
\item The fractional flow $f=f(S_i)\in C^2([0,1])$ is monotone increasing. 
\end{enumerate}
\label{ass:assumptions}
\end{assumption}
\subsection{The Vertically Integrated Model}
\label{sec:vertical-integrated}
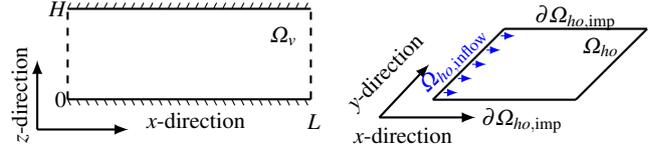
\begin{figure}
\hspace{-0.2cm}
\begin{tikzpicture}[scale = 0.4,>=latex]
\draw[->, thick] (-1,-01)--coordinate (x axis mid) (2,-01);
\draw[very thick] (2.3,-0.7) node[anchor=west] {$x$-direction}; 
\draw[->, thick] (-1,-1)--coordinate (y axis mid) (-1,1.3); 
\node[rotate=90, above=0cm] at (y axis mid) {$z$-direction};

\draw[thick] (-0.7,0) node[anchor=west]{$0$};
\draw[thick] (-0.9,3) node[anchor=west]{$H$};
\draw[thick] (8.1,-0.2) node[anchor=north]{$L$};
\draw[thick] (7.1,2.7) node[anchor=north]{$\Omega_v$};
\draw[thick] (0,0) -- (8,0);
\draw[thick] (0,3) -- (8,3); 

\draw[dashed, thick] (0,0) -- (0,3);
\draw[dashed, thick] (8,0) -- (8,3); 

\draw[-] (0,0) -- (0.1,-0.2);\draw[-] (0.3,0) -- (0.4,-0.2);\draw[-] (0.6,0) -- (0.7,-0.2);\draw[-] (0.9,0) -- (1,-0.2);\draw[-] (1.2,0) -- (1.3,-0.2);\draw[-] (1.5,0) -- (1.6,-0.2);\draw[-] (1.8,0) -- (1.9,-0.2);\draw[-] (2.1,0) -- (2.2,-0.2);\draw[-] (2.4,0) -- (2.5,-0.2);\draw[-] (2.7,0) -- (2.8,-0.2);\draw[-] (3,0) -- (3.1,-0.2);\draw[-] (3.3,0) -- (3.4,-0.2);\draw[-] (3.6,0) -- (3.7,-0.2);\draw[-] (3.9,0) -- (4,-0.2);\draw[-] (4.2,0) -- (4.3,-0.2);\draw[-] (4.5,0) -- (4.6,-0.2);\draw[-] (4.8,0) -- (4.9,-0.2);
\draw[-] (5.1,0) -- (5.2,-0.2);\draw[-] (5.4,0) -- (5.5,-0.2);\draw[-] (5.7,0) -- (5.8,-0.2);\draw[-] (6,0) -- (6.1,-0.2);\draw[-] (6.3,0) -- (6.4,-0.2);\draw[-] (6.6,0) -- (6.7,-0.2);\draw[-] (6.9,0) -- (7,-0.2);\draw[-] (7.2,0) -- (7.3,-0.2);\draw[-] (7.5,0) -- (7.6,-0.2);\draw[-] (7.8,0) -- (7.9,-0.2);

\draw[-] (0,3) -- (0.1,3.2);\draw[-] (0.3,3) -- (0.4,3.2);\draw[-] (0.6,3) -- (0.7,3.2);\draw[-] (0.9,3) -- (1,3.2);\draw[-] (1.2,3) -- (1.3,3.2);\draw[-] (1.5,3) -- (1.6,3.2);\draw[-] (1.8,3) -- (1.9,3.2);\draw[-] (2.1,3) -- (2.2,3.2);\draw[-] (2.4,3) -- (2.5,3.2);\draw[-] (2.7,3) -- (2.8,3.2);\draw[-] (3,3) -- (3.1,3.2);\draw[-] (3.3,3) -- (3.4,3.2);\draw[-] (3.6,3) -- (3.7,3.2);\draw[-] (3.9,3) -- (4,3.2);\draw[-] (4.2,3) -- (4.3,3.2);\draw[-] (4.5,3) -- (4.6,3.2);\draw[-] (4.8,3) -- (4.9,3.2);\draw[-] (5.1,3) -- (5.2,3.2);\draw[-] (5.4,3) -- (5.5,3.2);\draw[-] (5.7,3) -- (5.8,3.2);\draw[-] (6,3) -- (6.1,3.2);\draw[-] (6.3,3) -- (6.4,3.2);\draw[-] (6.6,3) -- (6.7,3.2);\draw[-] (6.9,3) -- (7,3.2);\draw[-] (7.2,3) -- (7.3,3.2);\draw[-] (7.5,3) -- (7.6,3.2);\draw[-] (7.8,3) -- (7.9,3.2);
\end{tikzpicture}
\begin{tikzpicture}[scale = 0.47,>=latex]
 \draw[->, thick] (-1.5,-0.5)--coordinate (x axis mid) (1.2,-0.5);
 \draw[very thick] (0.8,-1) node[anchor=east] {$x$-direction}; 
 \draw[->, thick] (-1.5,-0.5)--coordinate (y axis mid) (-0.1,1); 
 \node[rotate=45, above=0.1cm] at (y axis mid) {$y$-direction};
 \draw[blue, thick] (0.2,1.5) node[rotate=45, anchor=north]{$\Omega_{ho,\text{inflow}}$};
\draw[thick] (0,0) -- (4,0);
\draw[thick] (4,0) -- (6,2) ;
\draw[thick] (6,2) -- (2,2) ;
\draw[thick] (2,2) -- (0,0) ;
\draw[->,blue] (0.3,0.2) -- (0.7,0.2) ;
\draw[->,blue] (0.7,0.6) -- (1.1,0.6) ;
\draw[->,blue] (1.1,1) -- (1.5,1) ;
\draw[->,blue] (1.5,1.4) -- (1.9,1.4) ;
\draw[->,blue] (1.9,1.8) -- (2.3,1.8) ;
\draw[thick] (4.7,2) node[anchor=north]{$\Omega_{ho}$} ;
\draw[thick] (4,2.9) node[anchor=north]{$\partial\Omega_{ho,\text{imp}}$} ;
\draw[thick] (2.5,0) node[anchor=north]{$\partial\Omega_{ho,\text{imp}}$} ;
\end{tikzpicture}
\caption{A vertical-cross section of a 3D domain with bottom $(0,L)\times\{0\}$ and top $(0,L)\times\{H\}$ (left). A horizontal-cross section $\Omega_{ho}$ of a three-dimensional domain $\Omega$ (right).}
\label{fig:V-H cross-section}
\end{figure}
In the following, we recall the Vertically Integrated model (VI-model) for two-phase flow, see for example \cite{Gasda2009,Gasda2011,Guo2014}. For this, we define the vertically averaged quantities
\begin{align}
\begin{array}{rl}
\hat{\phi}(x,y)=& \displaystyle{\int_{0}^{H}} \phi (x,y,z)\,dz,\hspace{0.5cm} \hat{\textbf{v}}_{ho,\alpha} = \displaystyle{\int_{0}^{H}}  \textbf{v}_{ho,\alpha}\,dz,\vspace{5pt} \\
  \hat{S}_{\alpha}(x,y,t)=&\dfrac{1}{\hat{\phi}}\displaystyle{\int_{0}^{H}} \phi (x,y,z) S_{\alpha}(x,y,z,t)\,dz,\vspace{5pt} \\
 \widehat{\textbf{K}}(x,y)=&\displaystyle{\int_{0}^{H}}  \textbf{K}(x,y,z)\,dz,\hspace{0.5cm}\hat{\lambda}_{\alpha}= \widehat{\textbf{K}}^{-1}\displaystyle{\int_{0}^{H}}  \textbf{K} \lambda_{\alpha}\,dz.
\end{array}
\label{eq:integ-variables}
\end{align}
The vector $\textbf{v}_{ho,\alpha}=(u_{\alpha},v_{\alpha})^T$ is the velocity vector in the horizontal directions. Assuming that the fluids are in vertical equilibrium means that the phases' pressures $p_{\alpha},\,\alpha\in\{i,d\}$, are essentially hydrostatic \cite{Guo2014}, i.e.
\begin{align}
 p_{\alpha}(x,y,z,t) = \hat{p}_{\alpha}(x,y,t)+\rho_{\alpha}g z,
 \label{eq:ref-pressure}
\end{align}
where $\hat{p}_{\alpha}$ is the phase reference pressure at the bottom of the medium. Integrating the two-phase flow model \eqref{eq:two-phase} along the vertical coordinate from $0$ to $H$ (Figure \ref{fig:V-H cross-section}) and using the vertical equilibrium assumption \eqref{eq:ref-pressure} yields
\begin{align}
 \begin{array}{rl}
\hat{\phi} \partial_t \hat{S}_{\alpha}+ \nabla_{ho}\cdot \hat{\textbf{v}}_{ho,\alpha}&=0,\vspace{5pt}\\
 \hat{\textbf{v}}_{ho,\alpha}  &= - \hat{\lambda}_{\alpha} \widehat{\textbf{K}} \nabla_{ho}\hat{p}_{\alpha},\vspace{5pt}\\
 \nabla_{ho}\cdot\hat{\textbf{v}}_{ho} &=0
\end{array}
\label{eq:integratedsystem}
\end{align}
in $\Omega_{ho} \times(0,T)$, where $\Omega_{ho}$ is a horizontal cross-section of the three-dimensional domain $\Omega$, see Figure \ref{fig:V-H cross-section}, with vertically averaged properties like porosity and permeability. We call this model Vertically Integrated model (VI-model). In this model, $\hat{\textbf{v}}_{ho} =\hat{\textbf{v}}_{ho,i}+\hat{\textbf{v}}_{ho,d}$ is the averaged total velocity and $\nabla_{ho}=(\partial_x,\,\partial_y)^T$ is the gradient operator in the horizontal $x$- and $y$-directions. The VI-model has two unknowns: the averaged saturation $\hat{S}_{i}=\hat{S}_{i}(x,y,t)\in[0,1]$ and the averaged pressure $\hat{p}_{i}=\hat{p}_{i}(x,y,t)\in\mathbb{R}$. Saturation and pressure of the defending phase can be again computed using the closure relations $\hat{S}_i + \hat{S}_d = 1$ and $\hat{p}_d - \hat{p}_i = \widehat{p_c}$, where $\widehat{p_c}=p_c(\hat{S}_i)$, respectively. 

The initial and boundary conditions in \eqref{eq:TPF-IBC} are also vertically averaged for the VI-model and are given as
\begin{align}
\begin{array}{rll}
 \hat{S}_i(\cdot,\cdot,0)&=0, \quad\quad &\text{ in } \Omega_{ho}, \\
 \hat{S}_i&=\hat{S}_{\text{inflow}},\quad\quad &\text{ in } \partial\Omega_{ho,\text{inflow}}\times [0,T],\\
 \textbf{n}\cdot\hat{\textbf{v}}_{ho,\alpha} &= 0, \quad\quad &\text{ on }\partial \Omega_{ho,\text{imp}}\times [0,T].
\end{array}
\label{eq:VI-IBC}
\end{align}
In \eqref{eq:VI-IBC}, we have $\hat{S}_{\text{inflow}}\in[0,1]$, and $\partial\Omega_{ho,\text{inflow}}$ is the inflow boundary of the horizontal cross-section $\Omega_{ho}$ (Figure \ref{fig:V-H cross-section}).

\subsection{The Multiscale Model}
\label{sec:multiscale}
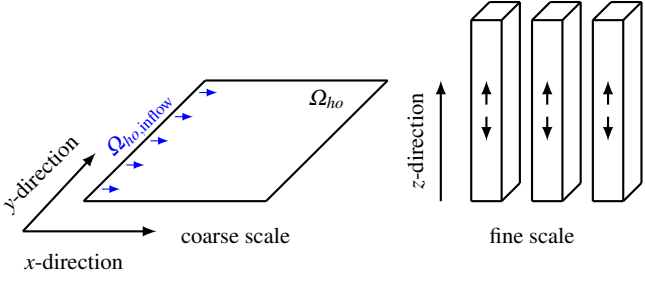
\begin{figure}
\centering
\begin{tikzpicture}[scale = 0.8,>=latex]
\hspace{-0.4cm}
  \draw[->, thick] (-1,-0.5)--coordinate (x axis mid) (1.2,-0.5);
 \draw[very thick] (0.8,-1) node[anchor=east] {$x$-direction}; 
 \draw[->, thick] (-1,-0.5)--coordinate (y axis mid) (0.2,0.8); 
 \node[rotate=45, above=0.1cm] at (y axis mid) {$y$-direction};
  \draw[blue, thick] (0.7,1.5) node[rotate=45, anchor=north]{$\Omega_{ho,\text{inflow}}$};

\draw[thick] (0,0) -- (3,0);
\draw[thick] (3,0) -- (5,2) ;
\draw[thick] (5,2) -- (2,2) ;
\draw[thick] (2,2) -- (0,0) ;
\draw[->,blue] (0.3,0.2) -- (0.6,0.2) ;
\draw[->,blue] (0.7,0.6) -- (1,0.6) ;
\draw[->,blue] (1.1,1) -- (1.4,1) ;
\draw[->,blue] (1.5,1.4) -- (1.8,1.4) ;
\draw[->,blue] (1.9,1.8) -- (2.2,1.8) ;
\draw[thick] (4,2) node[anchor=north]{$\Omega_{ho}$} ;
\draw[thick] (2.5,-0.3) node[anchor=north]{coarse scale} ;

\hspace{-0.5cm}
\draw[->, thick] (6.5,0)--coordinate (z axis mid) (6.5,2);
\node[rotate=90, above=0.1cm] at (z axis mid) {$z$-direction}; 
\draw[thick] (7,0) -- (7.5,0) ;
\draw[thick] (7.5,0) -- (7.8,0.3) ;
\draw[thick] (7,0) -- (7,3) ;
\draw[thick] (7.5,0) -- (7.5,3) ;
\draw[thick] (7.8,0.3) -- (7.8,3.3) ;
\draw[thick] (7,3) -- (7.5,3) ;
\draw[thick] (7.5,3) -- (7.8,3.3) ;
\draw[thick] (7,3) -- (7.3,3.3) ;
\draw[thick] (7.3,3.3) -- (7.8,3.3) ;
\draw[->, thick] (8.25,1.4) -- (8.25,1);
\draw[->, thick] (8.25,1.6) -- (8.25,2);

\draw[thick] (8,0) -- (8.5,0) ;
\draw[thick] (8.5,0) -- (8.8,0.3) ;
\draw[thick] (8,0) -- (8,3) ;
\draw[thick] (8.5,0) -- (8.5,3) ;
\draw[thick] (8.8,0.3) -- (8.8,3.3) ;
\draw[thick] (8,3) -- (8.5,3) ;
\draw[thick] (8.5,3) -- (8.8,3.3) ;
\draw[thick] (8,3) -- (8.3,3.3) ;
\draw[thick] (8.3,3.3) -- (8.8,3.3) ;
\draw[->, thick] (7.25,1.4) -- (7.25,1);
\draw[->, thick] (7.25,1.6) -- (7.25,2);

\draw[thick] (9,0) -- (9.5,0) ;
\draw[thick] (9.5,0) -- (9.8,0.3) ;
\draw[thick] (9,0) -- (9,3) ;
\draw[thick] (9.5,0) -- (9.5,3) ;
\draw[thick] (9.8,0.3) -- (9.8,3.3) ;
\draw[thick] (9,3) -- (9.5,3) ;
\draw[thick] (9.5,3) -- (9.8,3.3) ;
\draw[thick] (9,3) -- (9.3,3.3) ;
\draw[thick] (9.3,3.3) -- (9.8,3.3) ;
\draw[->, thick] (9.25,1.4) -- (9.25,1);
\draw[->, thick] (9.25,1.6) -- (9.25,2);
\draw[thick] (8,-0.3) node[anchor=north]{fine scale} ;
\end{tikzpicture}
\caption{An illustration of the two scales in the Multiscale model: the coarse-scale is the horizontal $xy$-plane while the fine-scale is a series of one-dimensional vertical domains. The arrows refer to the flow direction in the corresponding domains.  }
\label{fig:multiscale}
\end{figure}

Guo \emph{et al.} \cite{Guo2014} proposed a reduced model for asymptotically flat domains which they call the Multiscale model. It is used to describe $\text{CO}_2$ storage whenever the vertical equilibrium assumption \eqref{eq:ref-pressure} is not fully applicable. It comprises two scales: a coarse scale which is a horizontal cross-section $\Omega_{ho}$ of the three-dimensional domain $\Omega$, with vertically averaged properties, and a fine scale which is a series of one-dimensional vertical columns across the formation (Figure \ref{fig:multiscale}). The governing equations for this model are of two types: a two-dimensional equation for the vertically averaged pressure $\hat{p}_{i}$ (coarse-scale pressure \cite{Guo2014}), as in the VI-model \eqref{eq:integratedsystem}, and the two-phase flow model \eqref{eq:two-phase} for the saturation (fine-scale saturation \cite{Guo2014}), rearranged to focus on the vertical dynamics. These two type of equations are connected using an operator that reconstructs the coarse-scale pressure along the vertical direction. 

The Multiscale model is summarized for the invading phase $\alpha=i$ as follows:
\begin{align}
-\nabla_{ho}\cdot\widehat{\textbf{K}}\bigl( \lambda_{tot}(\hat{S}_i)  \nabla_{ho}\hat{p}_{i} \bigr)&=\nabla_{ho}\cdot\widehat{\textbf{K}}\bigl(\hat{\lambda}_{d}(\hat{S}_i)  \nabla_{ho}\widehat{p}_{c}(\hat{S}_i)\bigr)
\label{eq:multiscale1-a}
\end{align}
 in $\Omega_{ho}\times (0,T)$, and 
\begin{align}
\begin{array}{rl}
 p_{i}= & \hat{p}_{i}+ \pi_{i},\\\textbf{v}_{ho,i}=&-\lambda_{i}(S_i)\textbf{K}_{ho}\nabla p_{i},\\\partial_z w=&-\nabla_{ho}\cdot \textbf{v}_{ho} ,\\w_{i}= & f(S_i)w+f(S_i)\lambda_d(S_i)K_z\big(\partial_z p_c(S_i)\\& +(\rho_i-\rho_d)g\big),\\
 \phi\partial_{t}S_i+\nabla\cdot\textbf{v}_{i}=&0
\end{array}
\label{eq:multiscale1-b}
\end{align}
in $\Omega\times (0,T)$. Referencing to \cite{Guo2014}, a preferred choice of the reconstruction operator $\pi_{i}$, that is simple to compute and represents the natural forces in the medium, is given by 
\begin{align}
 \pi_i(x,y,z,t)=-&\displaystyle{\int_{0}^{z}} \Big( \big(S_i(x,y,r,t)\rho_i+ (1-S_i(x,y,r,t))\rho_d\big)g \nonumber\\+&S_i(x,y,r,t)\frac{\partial p_c\big(1-S_i(x,y,r,t)\big)}{\partial_z}\Big)dr,
 \label{eq:reconstruction}
\end{align}
for any $(x,y,z)\in\Omega$ and $t\in(0,T)$. The unknowns in this model are the vertically averaged pressure $\hat{p}_i=\hat{p}_i(x,y,t)\in \mathbb{R}$ and the saturation $S_i=S_i(x,y,z,t)\in[0,1]$. 

Note that the first equation in \eqref{eq:multiscale1-b} with the choice of $\pi_i$ in \eqref{eq:reconstruction} substitutes the vertical equilibrium assumption \eqref{eq:ref-pressure} in the VI-model. 
The Multiscale model \eqref{eq:multiscale1-a}, \eqref{eq:multiscale1-b} is closed with the initial and boundary conditions
\begin{align}
\begin{array}{rll}
 \hat{S}_i(\cdot,\cdot,0)&=0, \quad\quad &\text{ in } \Omega_{ho}, \\
 \hat{S}_i&=\hat{S}_{\text{inflow}},\quad\quad &\text{ in } \partial\Omega_{ho,\text{inflow}}\times [0,T],\\
  S_i(\cdot,\cdot,\cdot,0)&=0, \quad\quad &\text{ in } \Omega, \\
  S_i&=S_{\text{inflow}},\quad\quad &\text{ in } \partial\Omega_{\text{inflow}}\times [0,T],\\
 \textbf{n}\cdot\textbf{v}_i&= 0, \quad\quad &\text{ in }\partial \Omega_{\text{imp}}\times [0,T].
\end{array}
\label{eq:MS-IBC}
\end{align}

\subsection{The Vertical Equilibrium Model}
\label{sec:VE-model} 
Yortsos \cite{Yortsos} proposed the Vertical Equilibrium model (VE-model) to describe displacement processes in a vertical cross-section $\Omega_v=(0,L)\times(0,H)$ of an asymptotically flat domain $\Omega\subset \mathbb{R}^3$ such that the horizontal coordinate is parallel to the injection direction (Figure \ref{fig:displacementprocess}). 

We consider the model under the assumption of negligible gravity and capillary forces. First, we set $\gamma\coloneqq H/L$ and rescale the variables $x,y,t,K_x,K_z,p_{\alpha}$ in \eqref{eq:two-phase} into the dimensionless ones
\begin{equation}
  \begin{array}{rlrlrlrl}
\overline{x}=&\dfrac{x}{L},\quad&\overline{z}=&\dfrac{z}{H},\quad & \overline{t}=&\dfrac{t}{L/q}, \quad& \kappa_{j}=&\dfrac{K_{j}}{k_{j}}\\
\overline{u}_{\alpha}=&\dfrac{u_{\alpha}}{q},\quad &\overline{w}_{\alpha}=&\dfrac{w_{\alpha}}{q},\quad & \overline{p}_{\alpha}=&\dfrac{p_{\alpha}}{Lq\mu_{d}/k_{1}},  
   \end{array}
 \label{eq:dimensionlessvariables}
\end{equation}
for $j\in\{x,z\}$ and $\alpha\in\{i,d\}$. Here, $q>0$ is the inflow speed at the inflow boundary $\partial\Omega_{\text{inflow}}$ and $k_j$ is the mean value of the corresponding permeability function $K_j$. Then, after omitting the bar-signs, the TP-model \eqref{eq:two-phase} transforms into 
\begin{align}\left.
\begin{array}{rl}
 \phi\partial_{t}S_{\alpha} +\partial_{x}u_{\alpha}+\dfrac{1}{\gamma} \partial_{z}w_{\alpha} &=0,\\
 \partial_{x}u+\dfrac{1}{\gamma} \partial_{z}w &=0,\\
 u_{\alpha}&=-\lambda_{\alpha}(S_{\alpha})\, \kappa\, \partial_{x} p_{\alpha},\\
 \gamma w_{\alpha}&=-\lambda_{\alpha}(S_{\alpha})\, \kappa\, \partial_{ z}  p_{\alpha}
\end{array}\right.
\label{eq:Dimensionlesstwo-phase2}
\end{align}
in $D\times(0,T)$, where $D\coloneqq (0,1)\times(0,1)$ is the dimensionless spatial domain and $\kappa\coloneqq\kappa_x=\kappa_z$, by Assumption \ref{ass:assumptions}(1). Neglecting the capillary pressure implies $p_i=p_d=:p$ and one can show that the phases' velocities satisfy 
\begin{align*}
  u_{\alpha}=f(S_{\alpha})u,\quad\quad  w_{\alpha}=f(S_{\alpha})w.
\end{align*}
Then, the dimensionless TP-model \eqref{eq:Dimensionlesstwo-phase2} is reformulated in the fractional flow formulation with the unknown variables $S^{\gamma}_{\alpha},\,u^{\gamma},\,w^{\gamma}$ and $p^{\gamma}$ for some $\gamma>0$, 
\begin{align}\left.
\begin{array}{rl}
\phi \partial_{t}S^{\gamma}_{\alpha} + \partial_{x}(u^{\gamma}f(S^{\gamma}_{\alpha}))+\dfrac{1}{\gamma} \partial_{z}(w^{\gamma}f(S^{\gamma}_{\alpha})) &=0,\\
 \partial_{x}(u^{\gamma}) + \dfrac{1}{\gamma} \partial_{z}(w^{\gamma})&=0,\\
 u^{\gamma} &=-\lambda_{tot}(S^{\gamma}_{\alpha}) \kappa \partial_{x} p^{\gamma},\\
 \gamma w^{\gamma} &=-\lambda_{tot}(S^{\gamma}_{\alpha}) \kappa \partial_{z}p^{\gamma}
\end{array}\right.
\label{eq:Dimensionlessfractionalformulation}
\end{align}
in $D\times(0,T)$. Here, $u^{\gamma}$ and $w^{\gamma}$ are the velocity components in the horizontal- and vertical-directions. The unknowns in this model are pressure $p^{\gamma}$ and saturation $S^{\gamma}_{\alpha}$ for one of the phases. Yortsos applied formal asymptotic analysis for $\gamma\rightarrow 0$ and eliminated the pressure $p^{\gamma}$ from the velocity $(u^{\gamma},w^{\gamma})$ as clarified below.

We assume that for each $\gamma>0$ there exists a solution $(S^{\gamma},p^{\gamma},u^{\gamma},w^{\gamma})$ of \eqref{eq:Dimensionlessfractionalformulation}, with the asymptotic expansions
\begin{align}
 \begin{array}{ll}
  Z^{\gamma}=Z_{0}+\gamma Z_{1}+\mathcal{O}(\gamma^2),& \text{for }Z^{\gamma}\in\{S^{\gamma}_{\alpha},\,p^{\gamma},\,u^{\gamma}\},\vspace{3pt}\\
  w^{\gamma}=\gamma w_{1}+\mathcal{O}(\gamma^2).
 \end{array}
 \label{eq:asymptoticexpansion}
\end{align}
Note that the second equation in \eqref{eq:asymptoticexpansion} corresponds to the additional assumption, that the vertical velocity $w^{\gamma}$ is small. Using Assumptions \ref{ass:assumptions}(3) and \ref{ass:assumptions}(4), we have
\begin{align}
\begin{array}{rl}
 \lambda_{tot}(S^{\gamma}_{\alpha})=&\lambda_{tot}(S_{\alpha,0})+ \lambda_{tot}'(S_{\alpha,0})(\gamma S_{\alpha,1})+\mathcal{O}(\gamma^{2}),\vspace{3pt}\\ f(S^{\gamma}_{\alpha})=&f(S_{\alpha,0})+ f'(S_{\alpha,0})(\gamma S_{\alpha,1})+\mathcal{O}(\gamma^{2}).
\end{array}
 \label{eq:asymptoticexpansion1}
\end{align}
The second equation in \eqref{eq:Dimensionlessfractionalformulation} allows rewriting the first equation in \eqref{eq:Dimensionlessfractionalformulation} in a nonconservative form. Substituting the expansions \eqref{eq:asymptoticexpansion} and \eqref{eq:asymptoticexpansion1} into \eqref{eq:Dimensionlessfractionalformulation}, then considering terms of order $\mathcal{O}(1)$ in each equation leads to
\begin{align}
\begin{array}{rl}
\phi  \partial_{t}S_{\alpha,0} +u_{0}\partial_{x}f(S_{\alpha,0})+w_{1} \partial_{z}f(S_{\alpha,0})=&\mathcal{O}(\gamma),\vspace{3pt}\\
 \partial_{x} u_{0}+ \partial_{z} w_{1} =&\mathcal{O}(\gamma),\vspace{3pt}\\
-\lambda_{tot}(S_{\alpha,0}) \kappa \partial_{x} p_{0} =& u_{0},\vspace{3pt}\\
 -\lambda_{tot}(S_{\alpha,0}) \kappa \partial_{z} p_{0}=&\mathcal{O}(\gamma).
\end{array}
\label{eq:asym-system}
\end{align}
Due to the positivity of $\lambda_{tot}$ and $\kappa$ by Assumption \ref{ass:assumptions}(1) and \ref{ass:assumptions}(3), the last equation in \eqref{eq:asym-system} simplifies to
\begin{align}
  \partial_{z} p_{0}=\mathcal{O}(\gamma).
  \label{eq:derivedVE}
\end{align}
Thus, as $\gamma\rightarrow 0$, the global pressure $p_{0}$ is independent of the $z$-coordinate 
\begin{align}
 p_{0}=p_{0}(x,t).
 \label{eq:z-indep-pressure1}
\end{align}
The result in \eqref{eq:z-indep-pressure1} is a consequence of the asymptotic analysis and substitutes the vertical equilibrium assumption \eqref{eq:ref-pressure} in the VI-model.

In the following, the result \eqref{eq:z-indep-pressure1} is used to reformulate the velocity $\textbf{v}=(u,w)^T$ in terms of saturation only. Integrating the second equation in \eqref{eq:asym-system} along the vertical direction from $0$ to $1$, then using impermeability of the top and the bottom parts $\partial D_{\text{imp}}$ of the domain $\Omega$ yields
\begin{align}
 \partial_{x}\int_{0}^{1}u_{0}(x,z,t)\,dz=-\int_0^1 \partial_z w_1(x,z,t)\,dz=0,
\end{align}
as $\gamma\rightarrow 0$. Thus, we have
\begin{align}
 \int_{0}^{1}u_{0}dz=h,
 \label{eq:function-h}
\end{align}
for some function $h=h(t)$ with $h(t)>0$ for all $t\in[0,T]$. Integrating the third equation in \eqref{eq:asym-system} along the vertical direction from $0$ to $1$, then using equation \eqref{eq:function-h} gives
\begin{equation*}
-\int_{0}^{1}\lambda_{tot}(S_{\alpha,0}(x,z,t)) \kappa(x,z) \partial_{x} p_{0}(x,t)\,dz=h(t).
\end{equation*}
 Since $p_0$ is independent of the vertical coordinate $z$ we obtain that
\begin{equation}
\partial_{x} p_{0}(x,t)=-\dfrac{h(t)}{\int_{0}^{1}\lambda_{tot}(S_{\alpha,0}(x,z,t)) \kappa(x,z)\,dz}.
\label{eq:z-indep-pressure2}
\end{equation}
Substituting \eqref{eq:z-indep-pressure2} into the third equation in \eqref{eq:asym-system} gives a pressure independent formula for the horizontal velocity component $u_0$,
\begin{equation*}
 u_{0}[x,z,t;S_{\alpha,0}(x,z,t)]=\dfrac{h(t)\lambda_{tot}\bigl(S_{\alpha,0}(x,z,t)\bigr)\kappa(x,z)}{\int_{0}^{1}\lambda_{tot}\bigl(S_{\alpha,0}(x,r,t)\bigr) \kappa(x,r)\,dr},
\end{equation*}
for all $(x,z)\in D$ and $t\in(0,T)$. This formula is substituted in the second equation in \eqref{eq:asym-system} to produce a pressure independent formula for the vertical velocity component $w_1$,
\begin{equation*}
 w_{1}[x,z,t;S_{\alpha,0}(x,z,t)]=-h(t)\partial_{x}\dfrac{\int_{0}^{z}\lambda_{tot}\bigl(S_{\alpha,0}(x,r,t)\bigr)\kappa(x,r)\,dr}{\int_{0}^{1}\lambda_{tot}\bigl(S_{\alpha,0}(x,r,t)\bigr)\kappa(x,r)\,dr},
\end{equation*}
for all $(x,z)\in D$ and $t\in(0,T)$. Rescaling the time using $t\mapsto\bar{t}=\int_0^t h(r)dr+h(0)t$, then omitting the subscripts $\{0,1\}$ and the bar-signs yields the VE-model
\begin{align}
  \phi\partial_{t} S_{\alpha}+\partial_{x} \big(u[\cdot ,\cdot ;S_{\alpha}]&f(S_{\alpha})\big)+\partial_{z}\big(w[\cdot,\cdot;S_{\alpha}]~ f(S_{\alpha})\big)=0,\nonumber\\
  u[\cdot,\cdot;S_{\alpha}]=&\dfrac{\lambda_{tot}\bigl(S_{\alpha}\bigr)\kappa}{\int_{0}^{1}\lambda_{tot}\bigl(S_{\alpha}(\cdot,r,\cdot)\bigr)\kappa(\cdot,r)\,dr},\nonumber\\
  w[\cdot,z;S_{\alpha}(\cdot,z,\cdot)]=&-\partial_{x}\dfrac{\int_{0}^{z}\lambda_{tot}\bigl(S_{\alpha}(\cdot,r,\cdot)\bigr)\kappa(\cdot,r)\,dr}{\int_{0}^{1}\lambda_{tot}\bigl(S_{\alpha}(\cdot,r,\cdot)\bigr)\kappa(\cdot,r)\,dr}
 \label{eq:VEmodel}
\end{align}
in $D\times (0,T)$ for all $z\in (0,1)$. It is easy to check that the velocity vector $(u,w)^T$ in \eqref{eq:VEmodel} is divergence free. 

The velocity components in \eqref{eq:VEmodel} depend nonlocally on saturation. Hence, for the numerical comparisons in the next section, we suggest a finite-volume scheme with a non standard numerical flux that accounts for the nonlocality, and preserves the incompressibility property of the velocity.

\section{Comparison of the Darcy Models in Asymptotically Flat Domains}
\label{sec:VE-numericalresults}
In this section, we compare the VE-model \eqref{eq:VEmodel} with the other models from Section \ref{sec:background}. For this, we suggest a stable finite-volume scheme for the nonlocal VE-model. The comparisons show: firstly, the VE-model \eqref{eq:VEmodel} is more accurate than the VI-model, while being computationally more efficient than the full two-phase flow model. Secondly, the Multiscale model is equivalent to the VE-model in general.

\subsection{Finite-Volume Scheme for the VE-Model}
\label{sec:finitevolumescheme}
We consider a Cartesian grid for the dimensionless domain $D=(0,1)^2$,
	\begin{align*}
	 \mathcal{T}=\big\{T_{i,j}=&[x_{i-\frac{1}{2}},x_{i+\frac{1}{2}})\times[z_{j-\frac{1}{2}},z_{j+\frac{1}{2}})\mid(i,j)\in \mathcal{I}_x\times \mathcal{I}_z\subset \\ &  \mathbb{N}\times\mathbb{N}\big\},
	\end{align*}
	with the set of midpoints $\{(x_i,z_j)\mid(i,j)\in\mathcal{I}_x\times \mathcal{I}_z\}$. The set of edges of the cell $T_{i,j}$ is denoted by $\{E_{l}~|l\in \theta_{i,j}\}$ where
	 \begin{align*}
	  \theta_{i,j}:=\left\{\left(i-\frac{1}{2},j\right),\left(i+\frac{1}{2},j\right),\left(i,j-\frac{1}{2}\right),\left(i,j+\frac{1}{2}\right)\right\}.
	 \end{align*}
The set of neighbor cells of $T_{i,j}$ is defined as $\{T_{(i,j)_{l}}\mid l\in \theta_{i,j}\}$, where $T_{(i,j)_{l}}$ is the neighbor cell to $T_{i,j}$ with the common edge $E_{l}$. Let $N_x$ and $N_z$ be the number of cells in the horizontal and vertical directions, respectively. Then, the size of a cell is $|T_{i,j}|=\Delta x\Delta z$, where $\Delta x\coloneqq \tfrac{1}{N_x}$ and $\Delta z\coloneqq \tfrac{1}{N_z}$. For the initial condition $S_0\in[0,1]$, we define the cell-averaged values  
\begin{align}
 	 S_{i,j}^0 = \frac{1}{|T_{i,j}|}\int_{T_{i,j}} S_0(x,z)\,dx\,dz, \quad\quad (i,j)\in\mathcal{I}_x\times \mathcal{I}_z.
 	\end{align}
Given a positive integer $N$, the time interval $[0,T]$ is discretized into a set of disjoint subintervals $[t^{n},t^{n+1})$, with  $n\in\{0,1,2,...,N\}$, each of length $\Delta t>0$ such that $t^0=0$ and $t^N=T$. Then, a finite-volume scheme for the VE-model \eqref{eq:VEmodel} is given as
\begin{align}
& \frac{S_{i,j}^{n+1}-S_{i,j}^{n}}{\Delta t}+ \dfrac{1}{\Delta x \Delta z}\sum_{l\in\theta_{i,j}}\mathcal{F}_{l}\left(\bar{S}_{i,j}^{n},S_{i,j}^{n},S_{(i,j)_l}^{n}\right)=0.
\label{eq:VE-numerical scheme}
\end{align}
The crucial step in a finite-volume scheme is the choice of the numerical flux. We propose the numerical flux function
\begin{align}
	 \mathcal{F}_{l}\Big(\bar{S}_{i,j}^{n}&,S_{i,j}^{n},S_{(i,j)_l}^{n}\Big)= |E_l|\Big(\max\{\textbf{n}_{l}\cdot \textbf{v}_{l}[\bar{S}_{i,j}^{n}],0\}f(S_{i,j}^{n})\nonumber\\&+\min\{\textbf{n}_{l}\cdot \textbf{v}_{l}[\bar{S}_{i,j}^{n}],0\}f(S_{(i,j)_l}^{n})\Big),
	 \label{eq:numericalfluxfunction}
 	\end{align}
 	where $\textbf{n}_{l}$ is the outer normal to the edge $E_{l}$ of the cell $T_{i,j}$, $|E_l|$ is the length of the edge $E_l$ and $\bar{S}_{i,j}^n\coloneqq \{S_{r,m} \mid r\in\{i-1,i,i+1\},\,m\in\mathcal{I}_z\}$. The discrete velocity in \eqref{eq:numericalfluxfunction} $\textbf{v}_{l}^{n}[\bar{S}_{i,j}^{n}]=(u_{l}^{n}[\bar{S}_{i,j}^{n}],w_{l}^{n}[\bar{S}_{i,j}^{n}])^{T}$ is in fact nonlocal. Its components are defined as
 	\begin{align}
 	 u_{i+\frac{1}{2},j}^n=\,&\frac{1}{2}\Big(\frac{\lambda_{tot}(S_{i+1,j}^n)\kappa_{i+1,j}}{\Delta z\sum_{m=1}^{N_z} \lambda_{tot}(S_{i+1,m}^n)\kappa_{i+1,m}}\nonumber\\ & +\frac{\lambda_{tot}(S_{i,j}^n)\kappa_{i,j}}{\Delta z\sum_{m=1}^{N_z} \lambda_{tot}(S_{i,m}^n)\kappa_{i,m}}\Big)\quad\text{at }E_{i+\frac{1}{2},j},
 	 \label{eq:discretehorizontalvelocity1}
 	 \end{align}
 	 \begin{align}
 	 w_{i,j+\frac{1}{2}}^n=\,&-\frac{\Delta z}{\Delta x} \sum_{m=1}^j \left(u_{i+\frac{1}{2},m}^n-u_{i-\frac{1}{2},m}^n \right)\quad\text{at }E_{i,j+\frac{1}{2}}.
 	 \label{eq:discretehorizontalvelocity2}
 	\end{align}	
Due to the special discretization of the nonlocal terms \eqref{eq:discretehorizontalvelocity1} and \eqref{eq:discretehorizontalvelocity2} the scheme satisfies the following two important stability criteria:

  	\textbf{Mass conservation}:
	 Thanks to the definition of the numerical flux in \eqref{eq:numericalfluxfunction}, the finite-volume scheme \eqref{eq:VE-numerical scheme} is mass-conservative, i.e., if $E_l=E_{l'}$ such that $l\in \theta_{i,j}$ and $l'\in \theta_{(i,j)_l}$, then, for all $P_{i,j}\in[0,1]$, it holds that
	\begin{equation}
	  \mathcal{F}_l(\bar{P}_{i,j}, P_{i,j},P_{(i,j)_l})=-\mathcal{F}_{l'}(\bar{P}_{(i,j)_l},P_{(i,j)_l},P_{i,j}).
	\end{equation}

	\textbf{Discrete incompressibility}:
	If the vector of discrete velocity $\textbf{v}_{l}^{n}[\bar{S}_{i,j}^{n}]=(u_{l}^{n}[\bar{S}_{i,j}^{n}],~w_{l}^{n}[\bar{S}_{i,j}^{n}])^{T}$, $l\in\theta_{i,j}$, is defined as in \eqref{eq:discretehorizontalvelocity1} and \eqref{eq:discretehorizontalvelocity2}, then \emph{$\sum_{l\in\theta_{i,j}}\textbf{n}_{l}\cdot \textbf{v}_{l}^{n}[\bar{S}_{i,j}^{n}]=0$}.
	
\subsection{VE-Model vs. TP-Model}
\label{sec:VE-full}
In this section, we present several numerical examples to validate the VE-model and compare it with the TP-model \eqref{eq:Dimensionlessfractionalformulation} in domains with a decreasing geometrical parameter $\gamma\in\{1,\frac{1}{4},\frac{1}{8},\frac{1}{16},\frac{1}{32}\}$. The examples show also a remarkable reduction in computational complexity for the VE-model. For the comparisons, we consider the TP-model \eqref{eq:Dimensionlessfractionalformulation} with constant porosity $\phi=1$ and constant horizontal permeability $\kappa=1$. We apply the time-constant inflow condition
\begin{equation*}
S_{\text{inflow}}(z)=\left\{ 
\begin{array}{c l l}
	 0  &: & z\leq \frac{2}{5} \text{ and } z>\frac{3}{5},\\
	0.9 &: & \frac{2}{5}<z\leq \frac{3}{5}.
\end{array} \right.
\end{equation*}
The end time $T=0.3$ is chosen such that $S$ on the outflow boundary $\partial D_{\text{outflow}}$ vanishes in $[0,T]$. We set for the viscosity ratio $M=\mu_{d}/\mu_{i}=5$. Note that all numerical experiments are performed on the same hardware: (Intel(R) Core(TM) i7-4770 CPU @ 3.40GHz with memory 16GB).

\begin{figure}
\subfigure[TP-model $\gamma=1$]{
\includegraphics[scale=0.151]{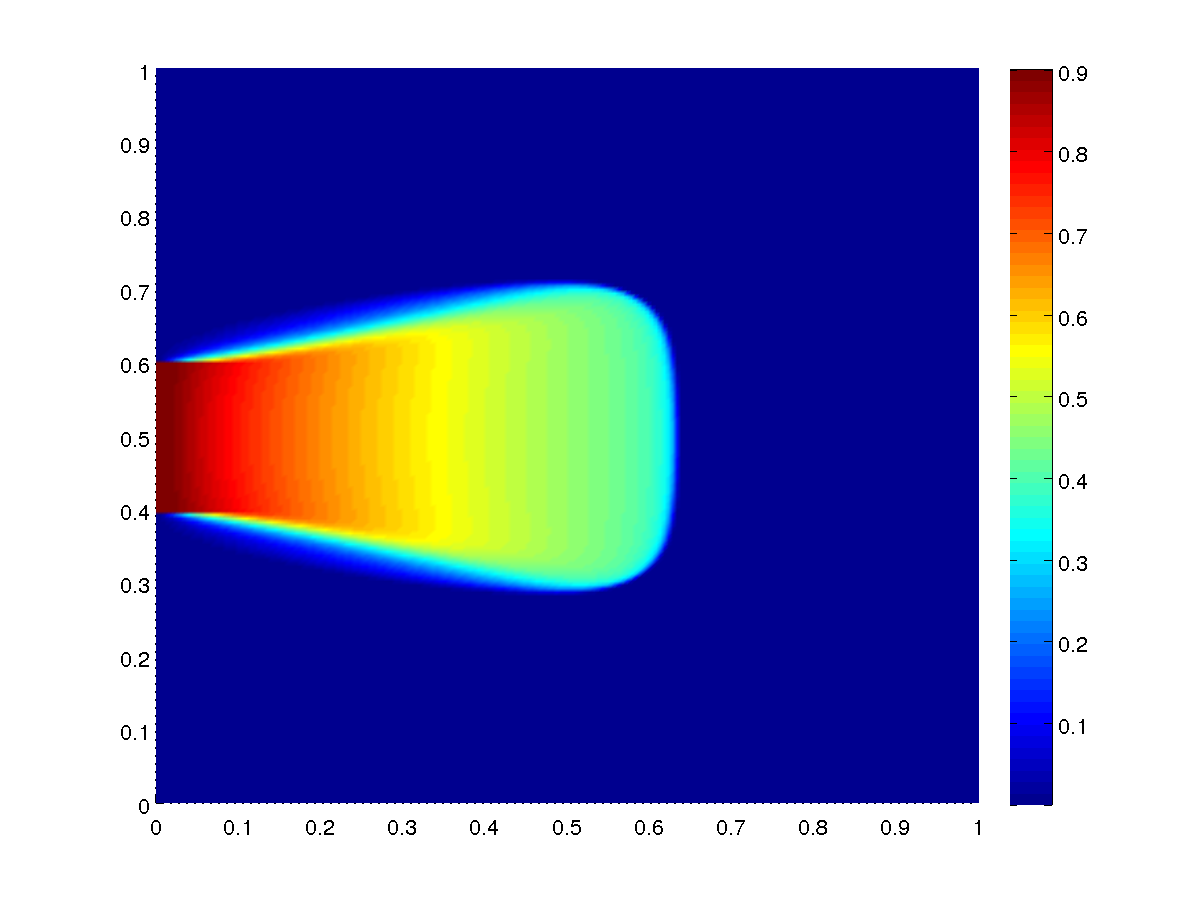}
\label{subfig:TP1}
}\hspace{-0.7cm}
\subfigure[TP-model $\gamma=\tfrac{1}{4}$]{
\includegraphics[scale=0.151]{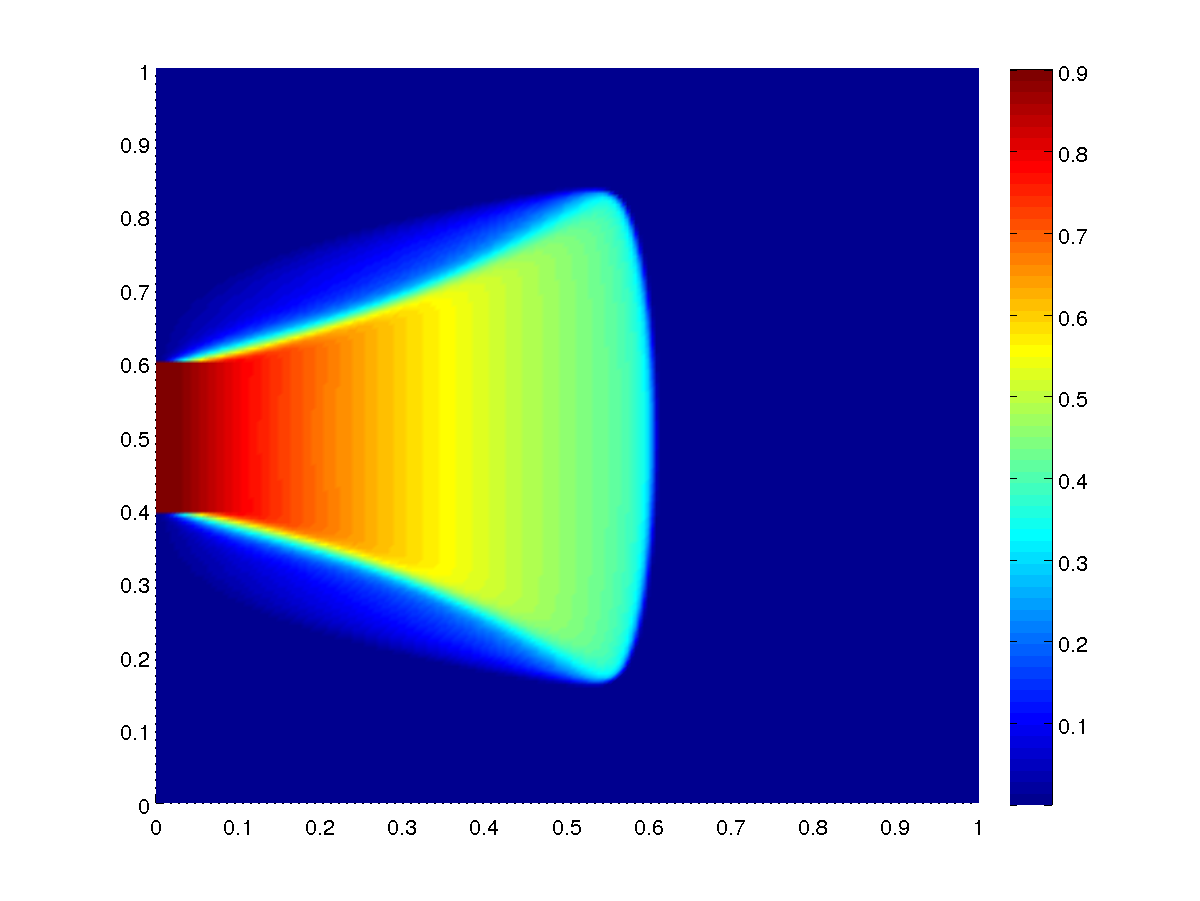}
\label{subfig:TP4}
}\hspace{-0.7cm}
\subfigure[TP-model $\gamma=\tfrac{1}{8}$]{
\includegraphics[scale=0.151]{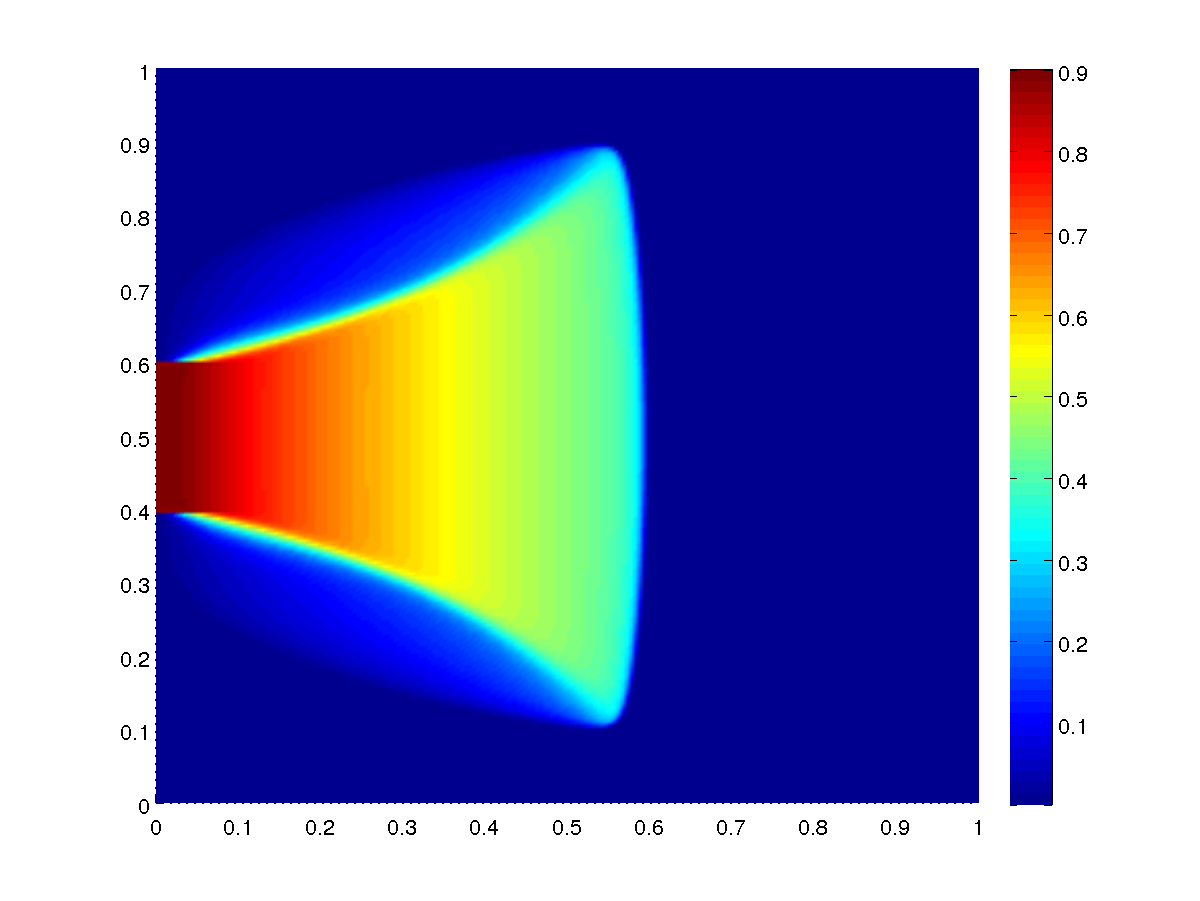}
\label{subfig:TP8}
}
\hspace{-0.5cm}
\subfigure[TP-model $\gamma=\tfrac{1}{16}$]{
\includegraphics[scale=0.151]{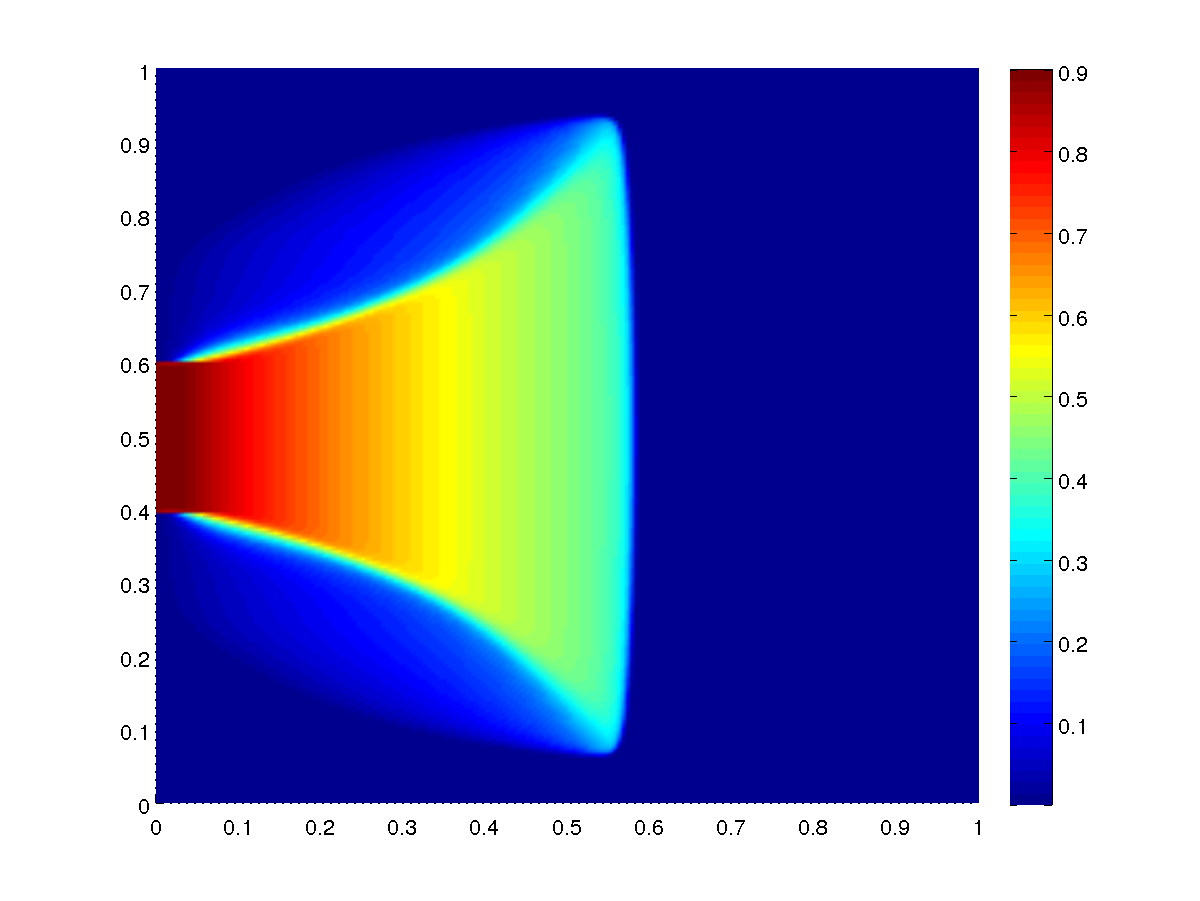}
\label{subfig:TP16}
}\hspace{-0.7cm}
\subfigure[TP-model $\gamma=\tfrac{1}{32}$]{
\includegraphics[scale=0.151]{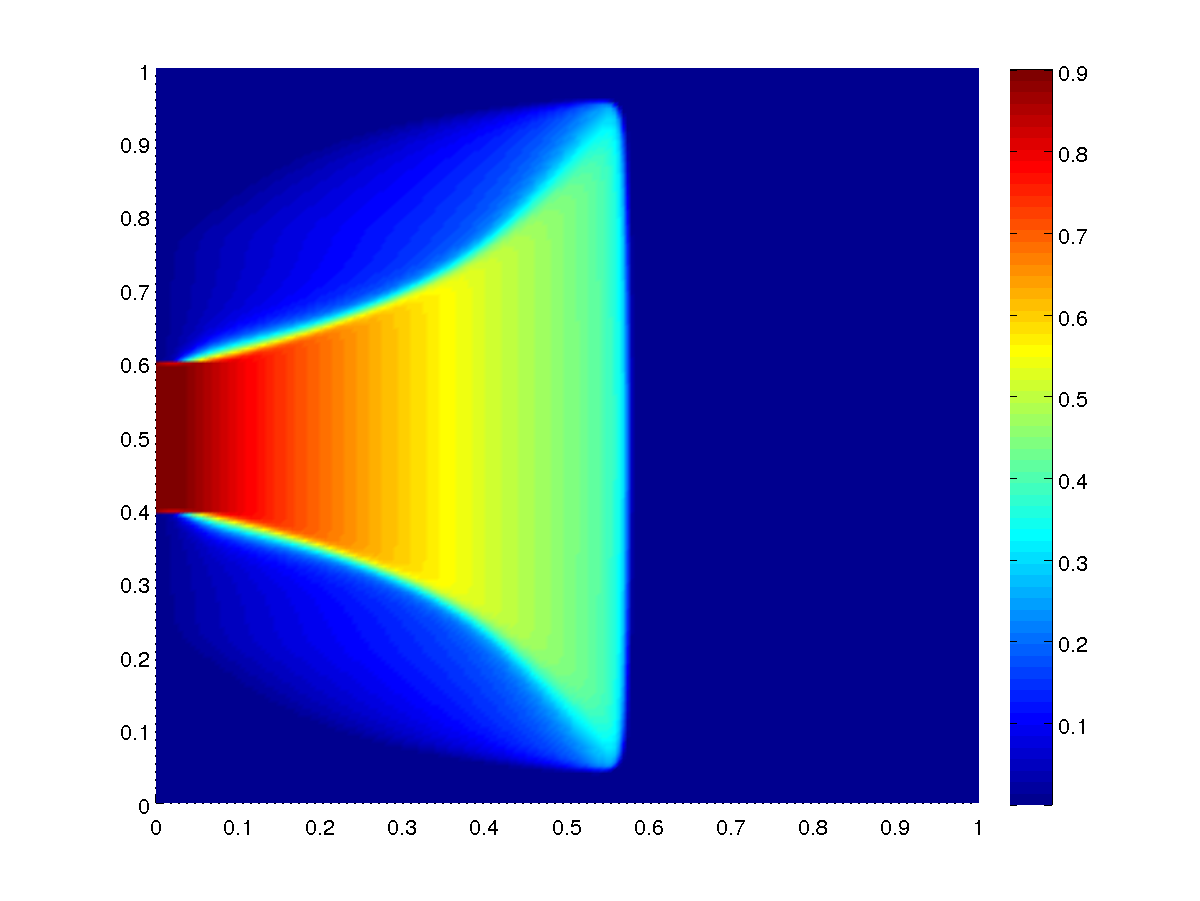}
\label{subfig:TP32}
}\hspace{-0.7cm}
\subfigure[VE-model]{
\includegraphics[scale=0.151]{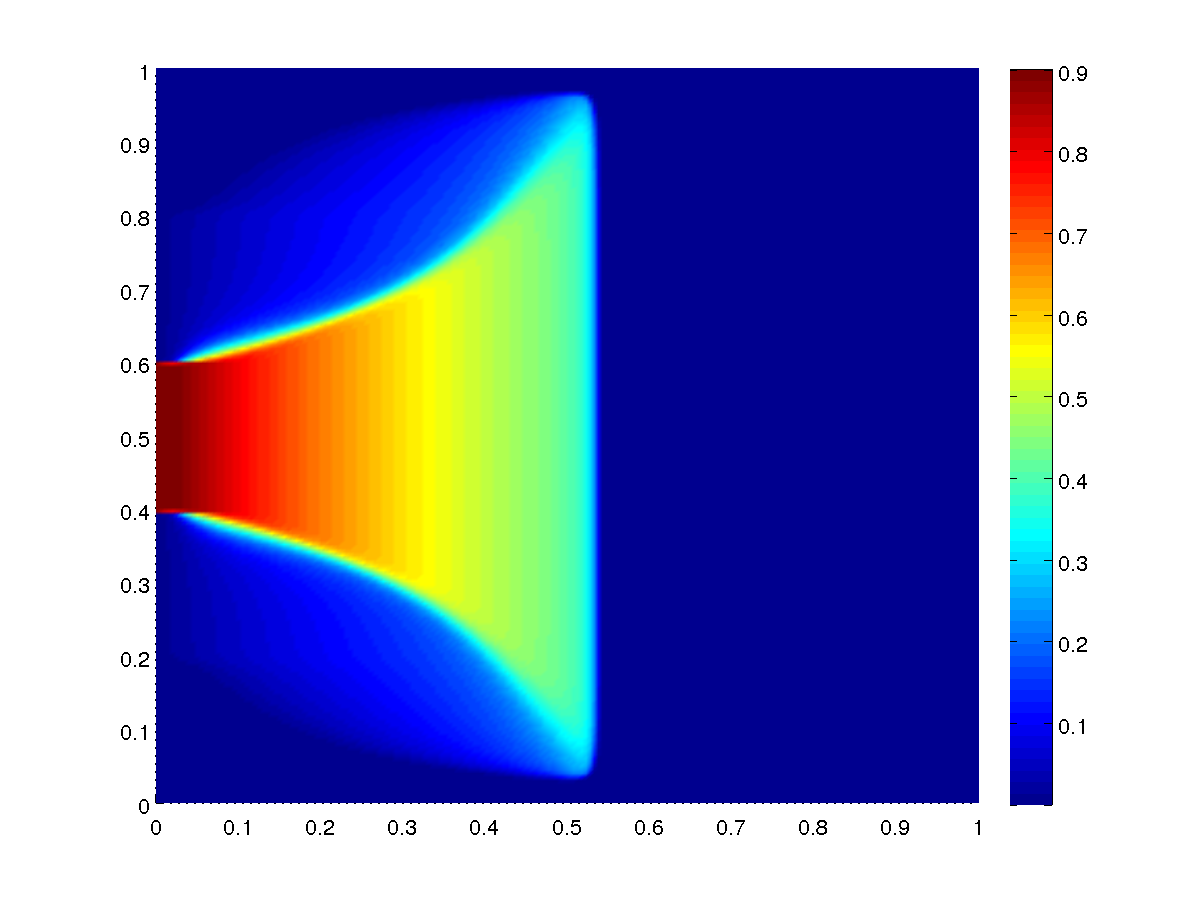}
\label{subfig:VE}
}
\caption{A comparison of the TP-model \eqref{eq:Dimensionlessfractionalformulation} for $\gamma\in\{1,\frac{1}{4},\frac{1}{8},\frac{1}{16},\frac{1}{32}\}$ with the VE-model \eqref{eq:VEmodel} using $M=5$ and $T=0.3$.}
\label{fig:numerics}
\end{figure}

\textbf{Example 1}: We consider the TP-model \eqref{eq:Dimensionlessfractionalformulation} in five domains $\Omega_{v,L}=(0,L)\times(0,1),~L\in\{1,4,8,16,32\}$. Figures \ref{subfig:TP1} - \ref{subfig:TP32} present the numerical solutions of the TP-model \eqref{eq:Dimensionlessfractionalformulation} using the IMPES-method (implicit pressure - explicit saturation) \cite{HuberHelmig1999} with pressure $p=1$ at the inflow boundary and $p=0$ at the outflow boundary. In Figure \ref{subfig:VE} the numerical solution of the VE-model \eqref{eq:VEmodel} is displayed, using the finite-volume scheme \eqref{eq:VE-numerical scheme}-\eqref{eq:discretehorizontalvelocity2}. The numerical solutions in Figure \ref{fig:numerics} correspond to a uniform Cartesian grid of $200 \times 200$ elements in the dimensionless domain $D$.

\begin{table}
\centering
\begin{tabular}{|c |c |c |c |c |}
\hline
Grid size $N_x\times N_z$ & $100 \times 100$ & $200 \times 200$ & $400 \times 400$ \\ [0.5ex]
\hline
VE-model & 0.94 s & 10 s & 124.51 s  \\
TP-model & 6.7 s & 64.5 s & 722.79 s \\ [1ex]
\hline
\end{tabular}
\caption{CPU-time for solving the VE-model compared to the CPU-time for solving the TP-model}.
\label{table:VE-timecomp1}
\end{table}

Figure \ref{fig:numerics} shows that numerical solutions of the TP-model \eqref{eq:Dimensionlessfractionalformulation} converge to the corresponding numerical solution of the VE-model, as the domain parameter $\gamma=(1/L)$ tends to $0$. Moreover, the VE-model describes the vertical dynamics in the medium and, consequently, the spreading speed of the invading front is captured well. These observations indicate that the VE-model is a proper reduction of the two-phase flow model for asymptotically flat domains. 

\textbf{Example 2}: Let the dimensionless domain $D$ be discretized into a uniform Cartesian grid with equal numbers of vertical and horizontal cells $N_z=N_x$. Set $\gamma=1/60$ in the TP-model \eqref{eq:Dimensionlessfractionalformulation}, $T=0.3$ and $M=5$. In Table \ref{table:VE-timecomp1}, we present the CPU-time required by the TP-model model using the IMPES-method and the VE-model using the finite-volume scheme in Section \ref{sec:finitevolumescheme}. 

The data from Table \ref{table:VE-timecomp1} shows that the computational time required to solve the VE-model is almost six times less than for the TP-model \eqref{eq:Dimensionlessfractionalformulation}. This is a consequence of solving an algebraic system that results from discretizing the two-dimensional elliptic equation of pressure in the two-phase flow model \eqref{eq:Dimensionlessfractionalformulation}. The performance of the elliptic solver could of course be optimized. This would improve the results for the TP-model, but not change the overall result.

\textbf{Example 3}: The comparison in Example 2 is repeated with a fixed number of vertical cells, but varying the number of horizontal cells $N_x=iN_z,\,i\in\{0.5,\,1,\,2,\,4,\,8\}$. This scenario fits to the case of asymptotically flat domains.

From Table \ref{table:VE-timecomp2} we see that the computational time required for solving the VE-model is almost ten times less than for the TP-model, when $N_x =2\times N_z$. However, in the case $N_x=16\times N_z$, the computational time required for solving the VE-model is almost $18$ times less than for the TP-model.

\begin{table}
\centering
\begin{tabular}{|c| c| c|}
\hline
Grid size $N_x\times N_z$ & VE-model & TP-model\\
\hline
$50 \times 100$  & 0.35 s&1.5 s \\
$100 \times 100$& 0.94 s&6.7 s \\ 
$200 \times 100$ &  1.9 s &30.8 s\\
$400 \times 100$&7.3 s & 103 s  \\
$800\times 100$ &  25 s &441 s\\ [0.5ex]
\hline
\end{tabular}
\caption{CPU-time for solving the VE-model compared to that for solving the TP-model such that $N_x=i N_z,\,i\in \{0.5,\,1,\,2,\,4,\,8\}$.}
\label{table:VE-timecomp2}
\end{table}

\begin{figure}
\centering
 \includegraphics[scale=0.4]{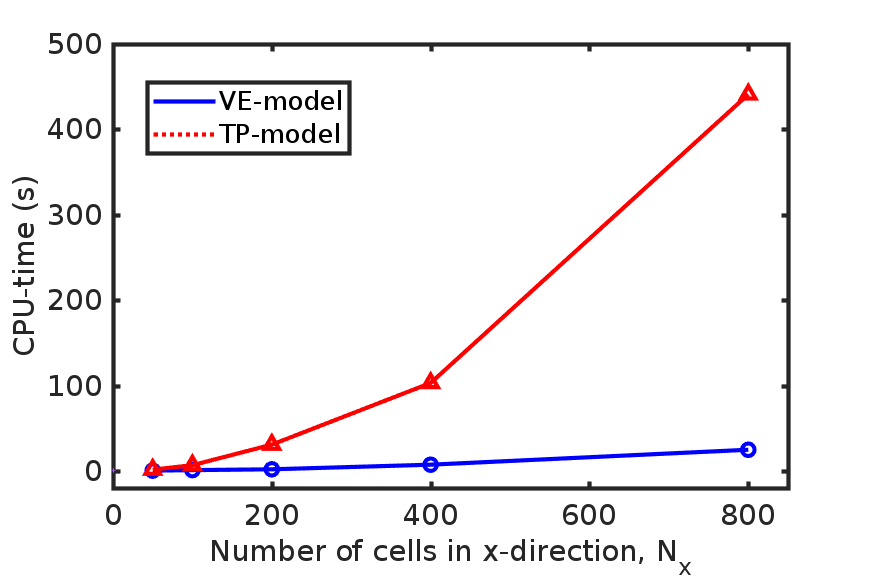}
 \caption{CPU-time for grids with fixed $N_z=100$ and varying $N_x\in\{50,\,100,\,200,\,400,\,800\}$.}
 \label{fig:CPU-time-VE-TPF}
\end{figure}

This example shows that increasing the number of horizontal cells $N_x$, while keeping $N_z$ fixed, limitedly increases the computational time for solving the VE-model as illustrated in Figure \ref{fig:CPU-time-VE-TPF}. On the contrary, the computational time for solving the TP-model significantly increases as the dimension of the algebraic system, resulting from applying the finite difference method to the elliptic equation of pressure, grows.

The numerical examples in this section show the validity of the VE-model in asymptotically flat domains. In addition, they demonstrate a significant reduction in computational complexity of the VE-model, mainly when using grids that are finer in the horizontal than in the vertical direction. 

\subsection{VE-Model vs. VI-Model}
\label{sec:VE-VI}
This section highlights the property of describing vertical dynamics by the VE-model and its effect on well-estimating the spreading speed. Different numerical examples that link the VE-model \eqref{eq:VEmodel} with the VI-model \eqref{eq:integratedsystem} are performed. 

\begin{figure}
\centering
\begin{tikzpicture}[scale = 0.7,>=latex]
\draw[->, thick] (-1,-0.7)--coordinate (x axis mid) (2,-0.7);
\draw[very thick] (2.3,-0.7) node[anchor=west] {$x$-direction};
\draw[->, thick] (-1,-0.7)--coordinate (y axis mid) (-1,1.3); 
\node[rotate=90, above=0.2cm] at (y axis mid) {$z$-direction};
\draw[thick] (0,0) -- (9,0);
\draw[dashed, thick] (0,0.5)--(9,0.5);
\draw[dashed, thick] (0,1)--(9,1);
\draw[dashed, thick] (0,1.5)--(9,1.5);
\draw[dashed, thick] (0,2)--(9,2);
\draw[dashed, thick] (0,2.5)--(9,2.5);
\draw[thick] (0,3) -- (9,3); 

\draw[dashed, thick] (0,0) -- (0,3);
\draw[dashed, thick] (9,0) -- (9,3); 

\draw[-] (0,0) -- (0.1,-0.2);\draw[-] (0.3,0) -- (0.4,-0.2);\draw[-] (0.6,0) -- (0.7,-0.2);\draw[-] (0.9,0) -- (1,-0.2);\draw[-] (1.2,0) -- (1.3,-0.2);\draw[-] (1.5,0) -- (1.6,-0.2);\draw[-] (1.8,0) -- (1.9,-0.2);\draw[-] (2.1,0) -- (2.2,-0.2);\draw[-] (2.4,0) -- (2.5,-0.2);\draw[-] (2.7,0) -- (2.8,-0.2);\draw[-] (3,0) -- (3.1,-0.2);\draw[-] (3.3,0) -- (3.4,-0.2);\draw[-] (3.6,0) -- (3.7,-0.2);\draw[-] (3.9,0) -- (4,-0.2);\draw[-] (4.2,0) -- (4.3,-0.2);\draw[-] (4.5,0) -- (4.6,-0.2);\draw[-] (4.8,0) -- (4.9,-0.2);
\draw[-] (5.1,0) -- (5.2,-0.2);\draw[-] (5.4,0) -- (5.5,-0.2);\draw[-] (5.7,0) -- (5.8,-0.2);\draw[-] (6,0) -- (6.1,-0.2);\draw[-] (6.3,0) -- (6.4,-0.2);\draw[-] (6.6,0) -- (6.7,-0.2);\draw[-] (6.9,0) -- (7,-0.2);\draw[-] (7.2,0) -- (7.3,-0.2);\draw[-] (7.5,0) -- (7.6,-0.2);\draw[-] (7.8,0) -- (7.9,-0.2);\draw[-] (8.1,0) -- (8.2,-0.2);\draw[-] (8.4,0) -- (8.5,-0.2);\draw[-] (8.7,0) -- (8.8,-0.2);\draw[-] (9,0) -- (9.1,-0.2);

\draw[-] (0,3) -- (0.1,3.2);\draw[-] (0.3,3) -- (0.4,3.2);\draw[-] (0.6,3) -- (0.7,3.2);\draw[-] (0.9,3) -- (1,3.2);\draw[-] (1.2,3) -- (1.3,3.2);\draw[-] (1.5,3) -- (1.6,3.2);\draw[-] (1.8,3) -- (1.9,3.2);\draw[-] (2.1,3) -- (2.2,3.2);\draw[-] (2.4,3) -- (2.5,3.2);\draw[-] (2.7,3) -- (2.8,3.2);\draw[-] (3,3) -- (3.1,3.2);\draw[-] (3.3,3) -- (3.4,3.2);\draw[-] (3.6,3) -- (3.7,3.2);\draw[-] (3.9,3) -- (4,3.2);\draw[-] (4.2,3) -- (4.3,3.2);\draw[-] (4.5,3) -- (4.6,3.2);\draw[-] (4.8,3) -- (4.9,3.2);\draw[-] (5.1,3) -- (5.2,3.2);\draw[-] (5.4,3) -- (5.5,3.2);\draw[-] (5.7,3) -- (5.8,3.2);\draw[-] (6,3) -- (6.1,3.2);\draw[-] (6.3,3) -- (6.4,3.2);\draw[-] (6.6,3) -- (6.7,3.2);\draw[-] (6.9,3) -- (7,3.2);\draw[-] (7.2,3) -- (7.3,3.2);\draw[-] (7.5,3) -- (7.6,3.2);\draw[-] (7.8,3) -- (7.9,3.2);\draw[-] (8.1,3) -- (8.2,3.2);\draw[-] (8.4,3) -- (8.5,3.2);\draw[-] (8.7,3) -- (8.8,3.2);\draw[-] (9,3) -- (9.1,3.2);

\end{tikzpicture}
\caption{A vertical-cross section of a 3D domain consisting of $6$ thin layers.}
\label{fig:j-layers}
\end{figure}
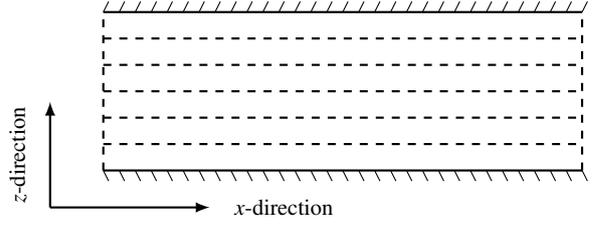

In a vertical cross-section $\Omega_v$ of $\Omega\subset \mathbb{R}^3$, the VI-model \eqref{eq:integratedsystem} reduces to a one-dimensional equation. Using dimensionless variables it is given by
\begin{align}
 \partial_t S+ \partial_z f(S)=0 \quad\quad\text{ in } (0,1)\times(0,T).
 \label{eq:VI}
\end{align}
We establish a comparison between the VI-model \eqref{eq:VI} and the VE-model \eqref{eq:VEmodel} such that the two-dimensional domain $D=(0,1)\times(0,1)$ consists of a set of adjacent horizontal layers as in Figure \ref{fig:j-layers}. This reduces the VE-model to a coupled system of one-dimensional equations, each of them describing the dynamics in the corresponding layer and taking into account the mass vertical exchange between adjacent layers. Then, the number of vertical cells $N_z$ in the finite-volume scheme \eqref{eq:VE-numerical scheme} equals the number of layers in the domain.

It is easy to check that the VI-model \eqref{eq:VI} corresponds to the one-layer scenario $N_z=1$ in the VE-model. Therefore, the comparison of the VE-model \eqref{eq:VEmodel} to the VI-model \eqref{eq:VI} is based on varying the number of layers in the domain: $N_z=1$ (the VI-model), $N_z=2$, and $N_z=5$. The two-dimensional VE-model is used here as a reference for accuracy. 

It is clear that the computational complexity of the VE-model with $N_z=1$ (the VI-model) is the lowest, and increases with the number of layers. However, the following examples show that decreasing the number of layers to $N_z=1$ has a strong impact on overestimating the spreading speed of the invading fluid. 

The models' comparison is based on two examples. In the first example, $S_{\text{inflow}}$ is independent of $z$, but the permeability $\kappa$ changes along the vertical axis. In the second example, $\kappa$ is $z$-independent but $S_{\text{inflow}}$ is injected only at the lower part of the inflow boundary. In both examples, we set $M=2$, $N_x=1000$ and $T=0.3$. 

\begin{figure}
\centering
\subfigure[VI model \eqref{eq:VI}]{
\includegraphics[scale=0.19]{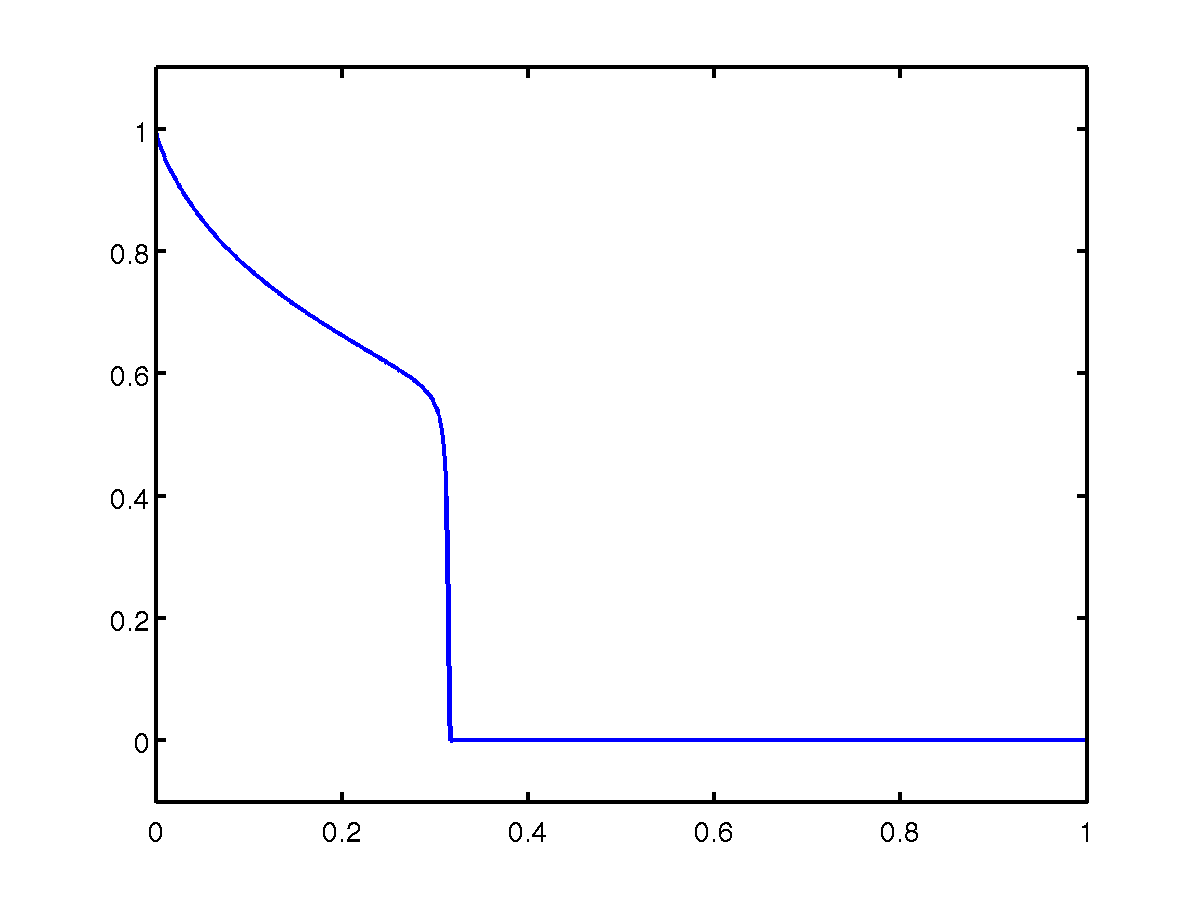}
\label{subfig:perm1-L1}
}
\subfigure[VE-model \eqref{eq:VEmodel} with $N_z=2$]{
\includegraphics[scale=0.19]{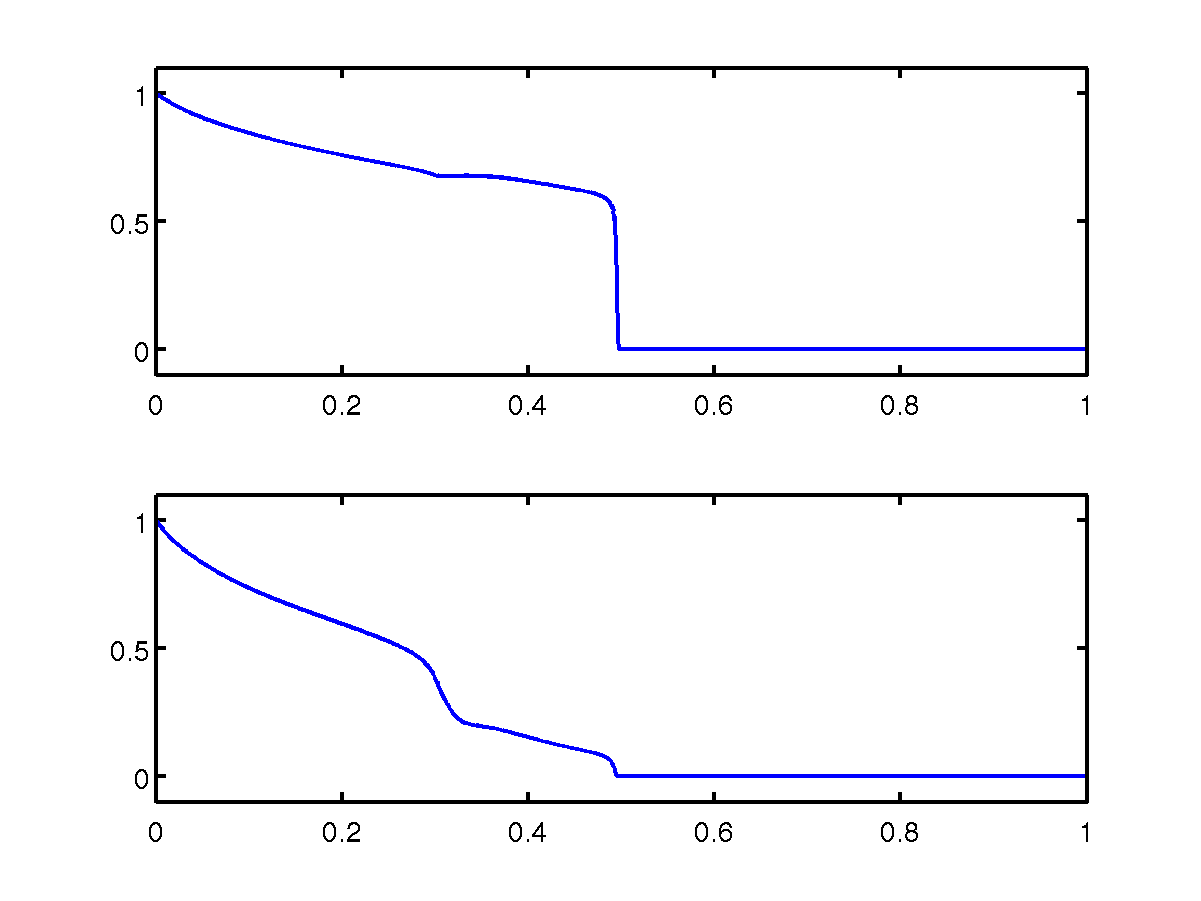}
\label{subfig:perm1-L2}
}\\
\subfigure[VE-model \eqref{eq:VEmodel} with $N_z=5$]{
\includegraphics[scale=0.19]{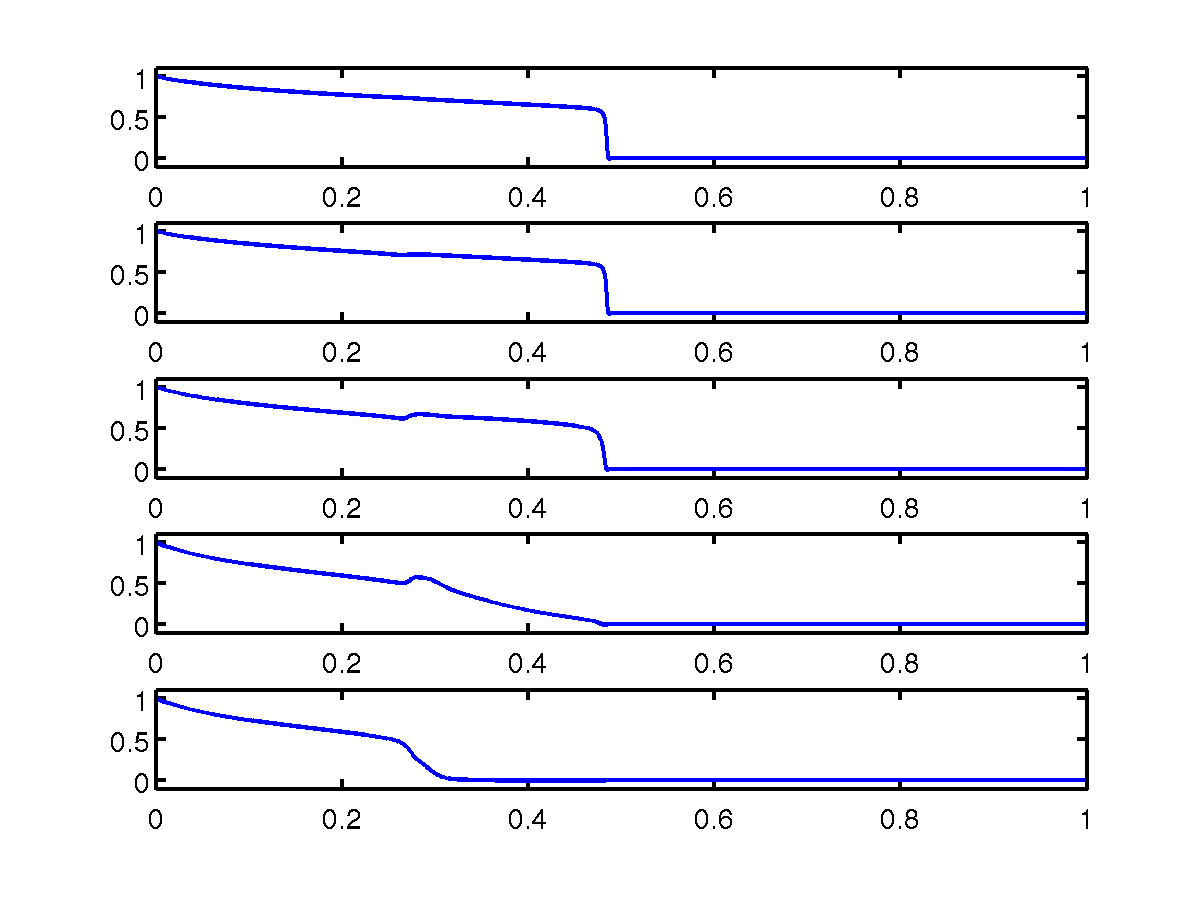}
\label{subfig:perm1-L5}
}
\subfigure[VE-model \eqref{eq:VEmodel} with $N_z=100$]{
\includegraphics[scale=0.19]{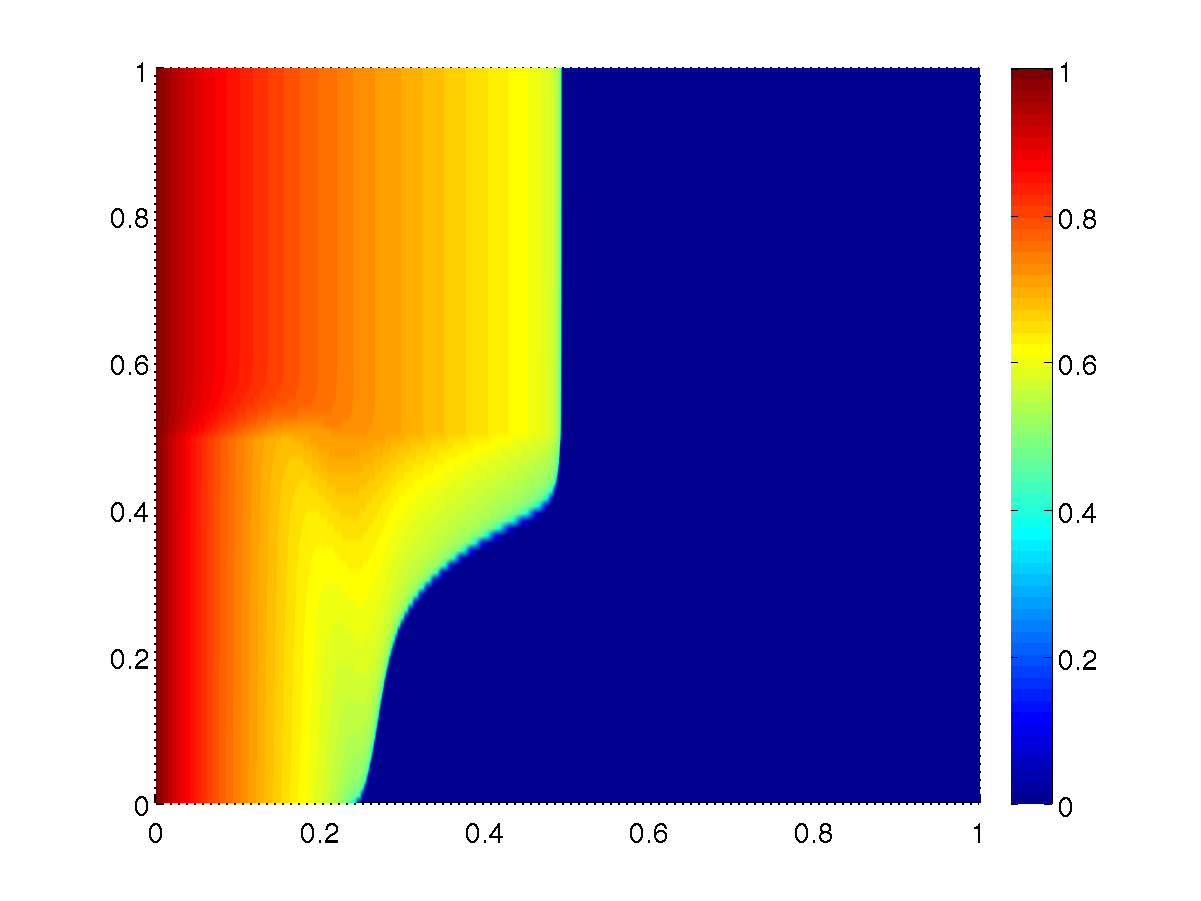}
\label{subfig:perm1-VE}
}
\caption{A comparison of the VI-model \eqref{eq:VI} to the VE-model \eqref{eq:VEmodel} for $\kappa$ as in \eqref{eq:perm1-2}, $S_{\text{inflow}}=1$, $M=2$, $T=0.3$ and $\Delta x=0.001$.}
\label{fig:perm1}
\end{figure}

\textbf{Example 1}:
We assume constant saturation $S_{\text{inflow}}=1$ along the inflow boundary. The domain consists of two parts with different permeabilities
\begin{align}
 \kappa=\left\{\begin{array}{cc}
                1,\quad & \text{in } D_{\text{upper}}:=(0,1)\times(0.5,1),\\
                0.5,\quad &  \text{in } D_{\text{lower}}:=(0,1)\times (0,0.5).
               \end{array}\right.
               \label{eq:perm1-2}
\end{align}
This condition is approximated according to the number of layers, i.e.
\begin{align}
 \kappa_{j}= \frac{1}{\Delta z} \int_{(j-1)\Delta z}^{j\Delta z}\kappa(z)\,dz, \quad\quad j=1,\cdots, N_z.
\end{align}
In Figure \ref{fig:perm1}, we show the effect of different permeabilities in the horizontal layers on the models' accuracy. Figure \ref{subfig:perm1-L1} presents the numerical solution of the VI-model \eqref{eq:VI} with averaged permeability $\kappa=0.75$. In Figure \ref{subfig:perm1-L2}, the domain consists of two layers: the lower layer has low permeability $\kappa_1=0.5$ and the upper layer has high permeability $\kappa_2=1$. In Figure \ref{subfig:perm1-L5}, the domain consists of five layers. The averaged permeabilities in the layers are: $\kappa_1=0.5,\,\kappa_2=0.5,\,\kappa_3=0.75,\,\kappa_4=1$ and $\kappa_5=1$, respectively from bottom to top. Figure \ref{subfig:perm1-VE} represents 
the numerical solution of the VE-model \eqref{eq:VEmodel} using vertical step size $\Delta z=0.01$. 

To sum up, Figure \ref{fig:perm1} shows that the spreading speed of the invading fluid is overestimated by the VI-model \eqref{eq:VI}, compared to the reference solution in Figure \ref{subfig:perm1-VE}. On the contrary, the spreading speed and the saturation distribution of the VE-model with $N_z=2$ and $N_z=5$ are well captured. This result is not limited to this example, but expected to hold for other cases where permeability changes along the vertical dimension of the medium.

\begin{figure}
\centering
\subfigure[VI-model \eqref{eq:VI}]{
\includegraphics[scale=0.19]{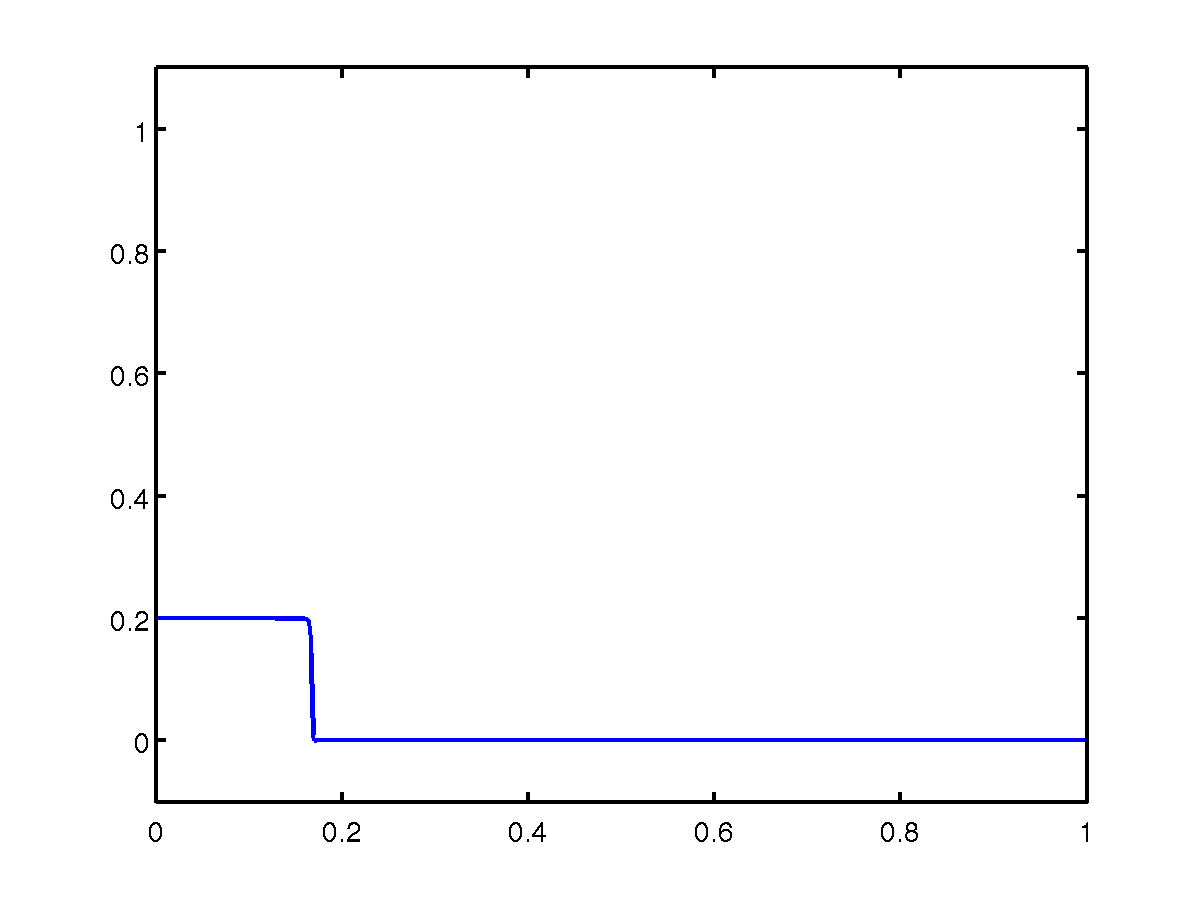}
\label{subfig:inflow-L1}
}
\subfigure[VE-model \eqref{eq:VEmodel} with $N=2$]{
\includegraphics[scale=0.19]{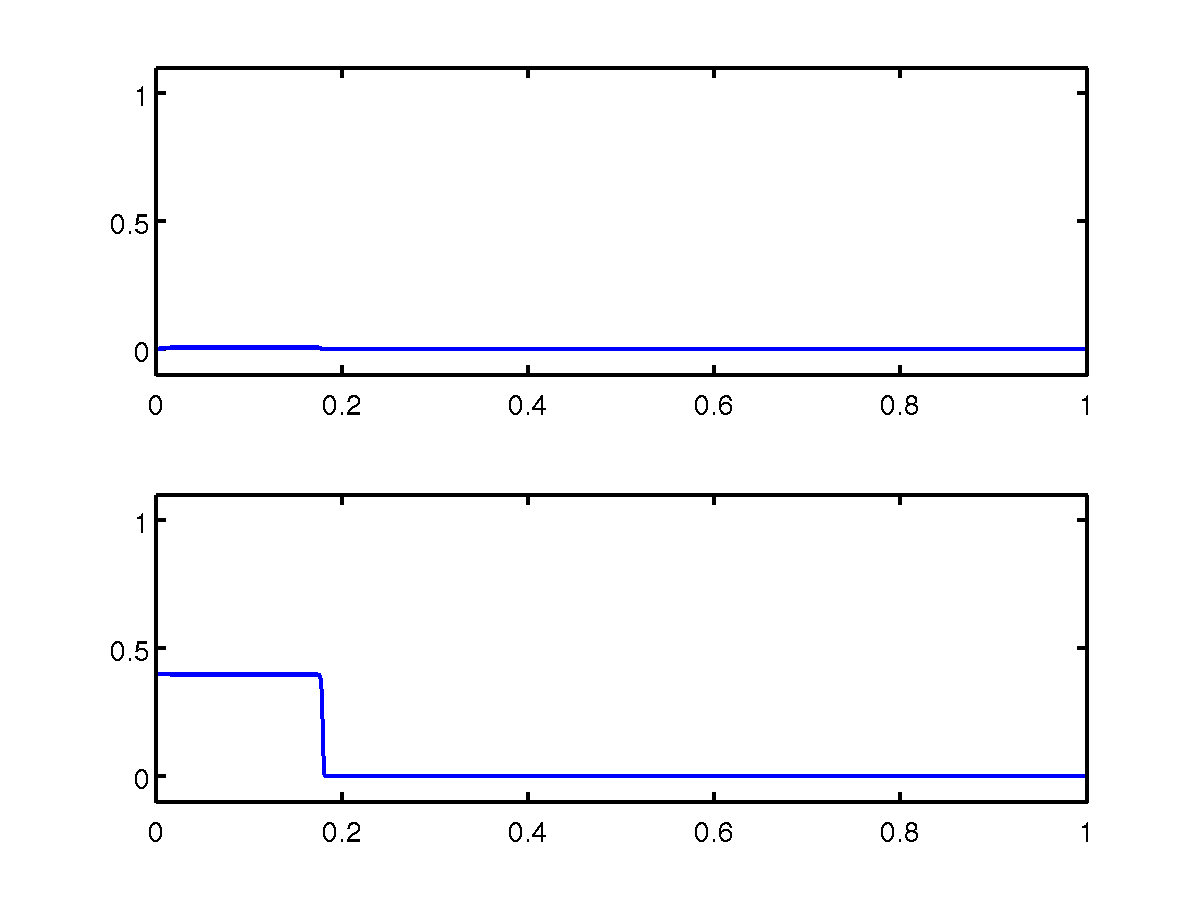}
\label{subfig:inflow-L2}
}
\subfigure[VE-model \eqref{eq:VEmodel} with $N=5$]{
\includegraphics[scale=0.19]{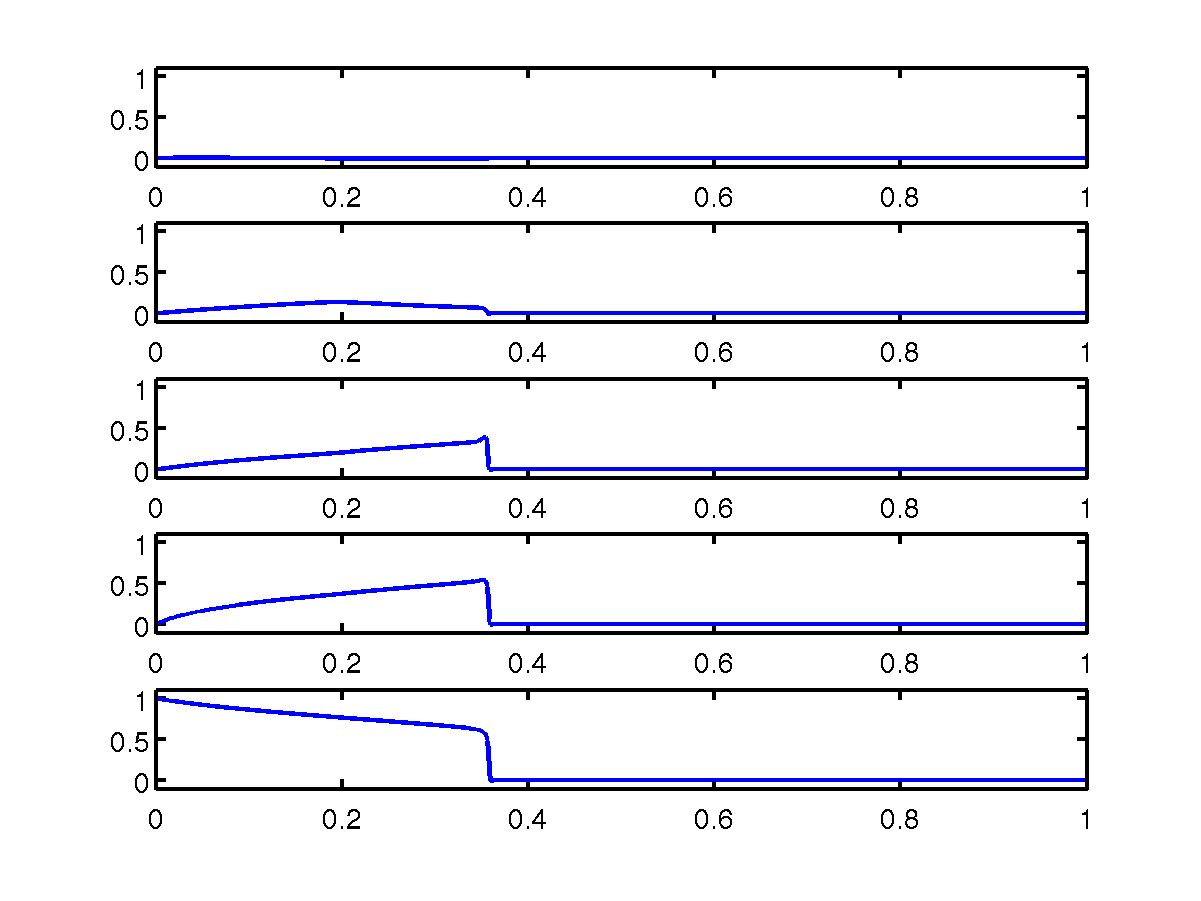}
\label{subfig:inflow-L5}
}
\subfigure[VE-model \eqref{eq:VEmodel} with $N=100$]{
\includegraphics[scale=0.19]{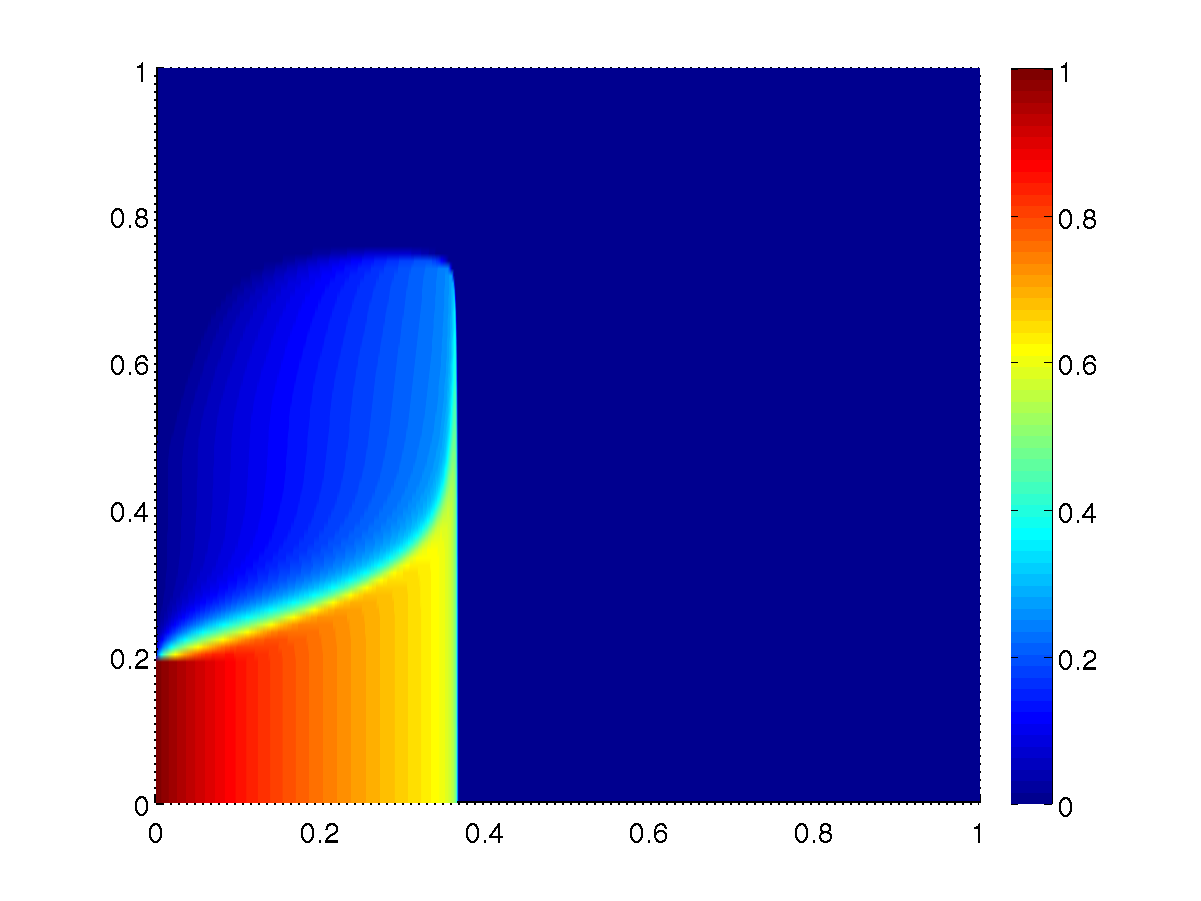}
\label{subfig:inflow-VE}
}
\caption{A comparison of the VI-model \eqref{eq:VI} with the VE-model \eqref{eq:VEmodel} for $S_{\text{inflow}}$ as in \eqref{eq:inflow1-2}, $\kappa=1$, $M=2$, $T=0.3$ and $\Delta x=0.001$.}
\label{fig:inflow}
\end{figure}

\textbf{Example 2}: We consider a domain with homogeneous permeability $\kappa=1$, such that the invading fluid is injected only at the lower part of the inflow boundary, i.e.
\begin{align}
 S_{\text{inflow}}(z)=\left\{\begin{array}{c c}
                          1, \quad& 0\leq z\leq 0.2,\\
                          0, \quad& 0.2\leq z\leq 1. 
                         \end{array}\right.
                         \label{eq:inflow1-2}
\end{align}
This condition is discretized according to the number of layers $N_z$ in the domain, i.e.
\begin{align}
 S_{\text{inflow},j}= \frac{1}{\Delta z} \int_{(j-1)\Delta z}^{j\Delta z} S_{\text{inflow}}(z)\,dz, \quad\quad j=1,\cdots, N_z.
\end{align}
In Figure \ref{subfig:inflow-L1}, we show the numerical solution of the VI-model \eqref{eq:VI} with averaged inflow $S_{\text{inflow}}=0.2$. Figure \ref{subfig:inflow-L2} provides the numerical solution of the VE-model with $N_z=2$ such that $S_{\text{inflow},1}=0.4$ in the lower layer and $S_{\text{inflow},2}=0$ in the upper one. The numerical solution of the VE-model with $N_z=5$ is presented in Figure \ref{subfig:inflow-L5} such that $S_{\text{inflow},1}=1,\,S_{\text{inflow},2}=0,\,S_{\text{inflow},3}=0,\,S_{\text{inflow},4}=0$ and $S_{\text{inflow},5}=0$ from bottom to top. Finally, the numerical solution of the VE-model \eqref{eq:VEmodel} with $\Delta z=0.01$ is given in Figure \ref{subfig:inflow-VE} as a reference for accuracy.  

This example emphasizes the weakness of the VI-model \eqref{eq:VI} whenever the fluid's injection position depends on the $z$-coordinate of the domain. In contrast to Example 2, the two-layers case $N_z=2$ in the VE-model is not sufficient to give a good estimation of the spreading speed. However, the spreading speed and the fluids' distribution in the five-layers case $N_z=5$ are in excellent agreement with the reference solution in Figure \ref{subfig:inflow-VE}.
\medskip 

Example 1 and Example 2 show the impact of varying the vertical permeability and the injection position on overestimating the spreading speed using the VI-model. On the contrary, the VE-model provides better estimates even using a very small number of vertical cells.

\subsection{VE-Model vs. Multiscale Model}
\label{sec:VE-Multi}
Under the same conditions assumed to derive the VE-model \eqref{eq:VEmodel}, we prove in this section that the Multiscale model \eqref{eq:multiscale1-a}, \eqref{eq:multiscale1-b} and \eqref{eq:reconstruction} in the vertical cross-section $\Omega_v=(0,1)\times(0,1)$ is equivalent to the VE-model. For the sake of simplicity, we also assume negligible gravity and capillary forces. Then, using the dimensionless variables \eqref{eq:dimensionlessvariables}, the Multiscale model \eqref{eq:multiscale1-a}, \eqref{eq:multiscale1-b} reduces to
\begin{align}
 -\partial_x \bigl( \hat{\lambda}_{tot}(\hat{S}_{\alpha})\hat{\kappa}\partial_x\hat{p}\bigr)=0 \quad\quad\text{ in }(0,1)\times(0,T),
 \label{eq:multiscale-dimensionless-a}
\end{align}
and in $D\times(0,T)$
\begin{align}
\begin{array}{rl}
  p(x,z,t)=&\,\hat{p}(x,t),\\
  -\lambda_{tot}(S_{\alpha})\kappa_1 \partial_x p=&u,\\
  \partial_z w=&-\partial_x u,\\
  \phi\partial_{t}S_{\alpha}+\partial_x \big(uf(S_{\alpha})\big)+\partial_z\big( wf(S_{\alpha})\big)=&\,0.
\end{array}
\label{eq:multiscale-dimensionless-b}
\end{align}
Note here that the assumption of negligible gravity and capillary forces implies that $\pi_{\alpha}=0$ (see \eqref{eq:reconstruction}). This model is closed with the initial and boundary conditions
\begin{align}
\begin{array}{rll}
 \hat{S}_{\alpha}&=\hat{S}_{0}, \quad\quad &\text{ in } (0,1), \\
 \hat{S}_{\alpha}&=\hat{S}_{\text{inflow}},\quad\quad &\text{ on } \{0\}\times [0,T],\\
  S_{\alpha}&=S_{0}, \quad\quad &\text{ in } D, \\
  S_{\alpha}&=S_{\text{inflow}},\quad\quad &\text{ on } \partial D_{\text{inflow}}\times [0,T],\\
 \textbf{n}\cdot\textbf{v}_{\alpha}&= 0, \quad &\text{ on }\partial D_{\text{imp}}\times [0,T].
\end{array}
\label{eq:MS-D-IBC}
\end{align}

\begin{theorem}
For $T>0$, let $(S_M,u_M,w_M,\hat{p}_M)$ be a smooth solution of the Multiscale model \eqref{eq:multiscale-dimensionless-a}, \eqref{eq:multiscale-dimensionless-b} in $D\times (0,T)$ with the initial and boundary conditions \eqref{eq:MS-D-IBC}. For $\overline{T}:=\int_0^T h(r)dr+ h(0)T$, where $h(t)=\int_0^1 u_M(0,z,t)\,dz$, let $(S_{VE},u_{VE},w_{VE})$ be a smooth solution of the VE-model \eqref{eq:VEmodel} in $D\times (0,\overline{T})$ with the initial and boundary conditions \eqref{eq:TPF-IBC}. Then,  $ S_M(\cdot,\cdot,t)=S_{VE}(\cdot,\cdot,\bar{t}),\, u_M(\cdot,\cdot,t)=u_{VE}(\cdot,\cdot,\bar{t})$ and $w_M(\cdot,\cdot,t)=w_{VE}(\cdot,\cdot,\bar{t})$ in $D$ where $\overline{t}:=\int_0^t h(r)dr+ h(0)t$ for any $t\in(0,T)$. Furthermore, there is a function $p_{VE}=p_{VE}(S_{VE})$ with
\begin{align*}
\partial_x p_{VE}(x,t)= \frac{1}{\int_{0}^{1}\kappa(x,z)\lambda_{tot}\big(S_{VE}(x,z,t)\big)\,dz}
\end{align*}
and $\hat{p}_{M}(\cdot,t)=p_{VE}(\cdot,\bar{t})$ in $(0,1)$. 
 \label{thm:VE-multi}
\end{theorem}

In view of the equivalence result of Theorem \ref{thm:VE-multi}, the finite-volume scheme in Section \ref{sec:finitevolumescheme} is also applicable for the Multiscale model \eqref{eq:multiscale-dimensionless-a}, \eqref{eq:multiscale-dimensionless-b}. This leads to a reduced computational complexity of the Multiscale model compared with that using a multiscale algorithm proposed in \cite{Guo2014}. This algorithm is analogous to the IMPES-method with an implicit treatment for the pressure at the coarse scale and an explicit treatment for the saturation at the fine scale.

\begin{proof}
Let $(S_M,u_M,w_M,\hat{p}_M)$ be a smooth solution of the Multiscale model \eqref{eq:multiscale-dimensionless-a}, \eqref{eq:multiscale-dimensionless-b}. Then, the coarse-scale equation \eqref{eq:multiscale-dimensionless-a} implies that the vertically averaged horizontal velocity,
 \begin{align}
 \hat{u}_M=-\hat{\lambda}_{tot}(\hat{S}_{M})\hat{\kappa}\partial_x\hat{p}_M,
  \label{eq:multi3}
\end{align}
satisfies 
\begin{align}
\partial_x \hat{u}_M=0.
 \label{eq:multi2}
\end{align}
Integrating \eqref{eq:multi2} over the $x$-direction from $0$ to an arbitrary $x\in (0,1)$ yields
\begin{align}
 \hat{u}_M(x,t)-h_M(t)=0,
 \label{eq:almost-relation1}
 \end{align}
where $h_M(t)=\hat{u}_M(0,t)$ is the averaged horizontal velocity at the inflow boundary $\partial D_{\text{inflow}}$ with $h_M(t)>0$ for all $t\in[0,T]$. Substituting \eqref{eq:almost-relation1} into \eqref{eq:multi3} and noting that the vertically-averaged quantities $\hat{\kappa}$ and $\hat{\lambda}_{tot}$ are defined as in \eqref{eq:integ-variables}, we obtain that 
\begin{align}
 \partial_x \hat{p}_M(x,t)=\frac{h_M(t)}{\int_{0}^{1}\kappa(x,z)\lambda_{tot}(S_{M}(x,z,t))\,dz},
\end{align}
which coincides with the definition of $\partial_x p_{VE}\coloneqq \partial_x p_{0}$ in equation \eqref{eq:z-indep-pressure2}. Using the first equation in the fine-scale model \eqref{eq:multiscale-dimensionless-b}, we obtain that
\begin{align}
 \partial_x p_M(x,z,t)=\frac{h_M(t)}{\int_{0}^{1}\kappa(x,z)\lambda_{tot}(S_{\alpha}(x,z,t))\,dz},
 \label{eq:almost-comp3}
\end{align}
which is in fact independent of the vertical coordinate $z$ as the reconstruction operator $\pi_i$ in this setting equals zero. Substituting equation \eqref{eq:almost-comp3} into the second equation in \eqref{eq:multiscale-dimensionless-b} yields the horizontal velocity $u_M$ as an expression of saturation alone
\begin{align}
 u_M[x,z,t;S_{M}]=-\frac{ h_M(t) \kappa(x,z)\lambda_{tot}\big(S_{M}(x,z,t)\big)}{\int_{0}^{1}\kappa(x,z)\lambda_{tot}\big(S_{M}(x,z,t)\big)\,dz}.
 \label{eq:almost-comp4}
\end{align}
Substituting \eqref{eq:almost-comp4} into the third equation in \eqref{eq:multiscale-dimensionless-b} then integrating over the vertical coordinate from $0$ to an arbitrary $z\in(0,1)$ yields the vertical velocity as an expression of saturation alone
\begin{align}
 w_M[x,z,t;S_{M}]=&-\partial_x \int_{0}^z u_M[x,r,t;S_{M}]\,dr.
 \label{eq:almost-comp5}
\end{align}
We substitute equation \eqref{eq:almost-comp4} and \eqref{eq:almost-comp5} into the fourth equation in the fine-scale model \eqref{eq:multiscale-dimensionless-b}, then rescale the time using  $t\mapsto \bar{t}=\int_0^t h_M(r)dr+h_M(0)t$. As a result, the Multiscale model simplifies to 
\begin{align}
  \partial_{t} S_M+\partial_{x} \big(u_M[\cdot,\cdot;&S_M]f(S_M)\big)+~\partial_{z}\big(w[\cdot,\cdot;S_M]~ f(S_M)\big)=0,\nonumber\\
  u_M[\cdot,\cdot;S_M]&=\dfrac{\lambda_{tot}(S_M)\kappa}{\int_{0}^{1}\lambda_{tot}\big(S_M(\cdot,z,\cdot)\big)\kappa(\cdot,z)\,dz},\nonumber\\
  w_M[\cdot,z;S_M(\cdot,z,\cdot)]&=-\partial_{x}\dfrac{\int_{0}^{z}\lambda_{tot}\big(S_M(\cdot,r,\cdot)\big)\kappa(\cdot,r)\,dr}{\int_{0}^{1}\lambda_{tot}\big(S_M(\cdot,r,\cdot)\big)\kappa(\cdot,r)\,dr},
  \label{eq:MSmodel-asVE}
\end{align}
for all $z\in(0,1)$, in $D\times(0,\overline{T})$. Note that this scaled time interval coincides with that in the VE-model as the function $h_M$ is equivalent to $h_{VE}$, defined as in \eqref{eq:function-h}.

The simplified version of the Multiscale model \eqref{eq:MSmodel-asVE} requires initial and boundary conditions with respect to the unknown $S_M$ only. Hence, it is closed with the initial and boundary conditions
\begin{align}
\begin{array}{rll}
  S_M&=S_{0}, \quad\quad &\text{ in } D, \\
  S_M&=S_{\text{inflow}},\quad\quad &\text{ on } \partial D_{\text{inflow}}\times [0,T],\\
 \textbf{n}\cdot\textbf{v}_M&= 0, \quad\quad &\text{ on }\partial D_{\text{imp}}\times [0,T].
\end{array}
\end{align}
Hence, a smooth solution $(S_M,u_M,w_M,\hat{p}_M)$ of the Multiscale model \eqref{eq:multiscale-dimensionless-a}, \eqref{eq:multiscale-dimensionless-b} is also a solution of the VE-model. Proving the other direction follows the same argumentation in reversed order.
\qed
\end{proof}

\begin{figure}
\centering
\subfigure[multiscale algorithm]{
\includegraphics[scale=0.2]{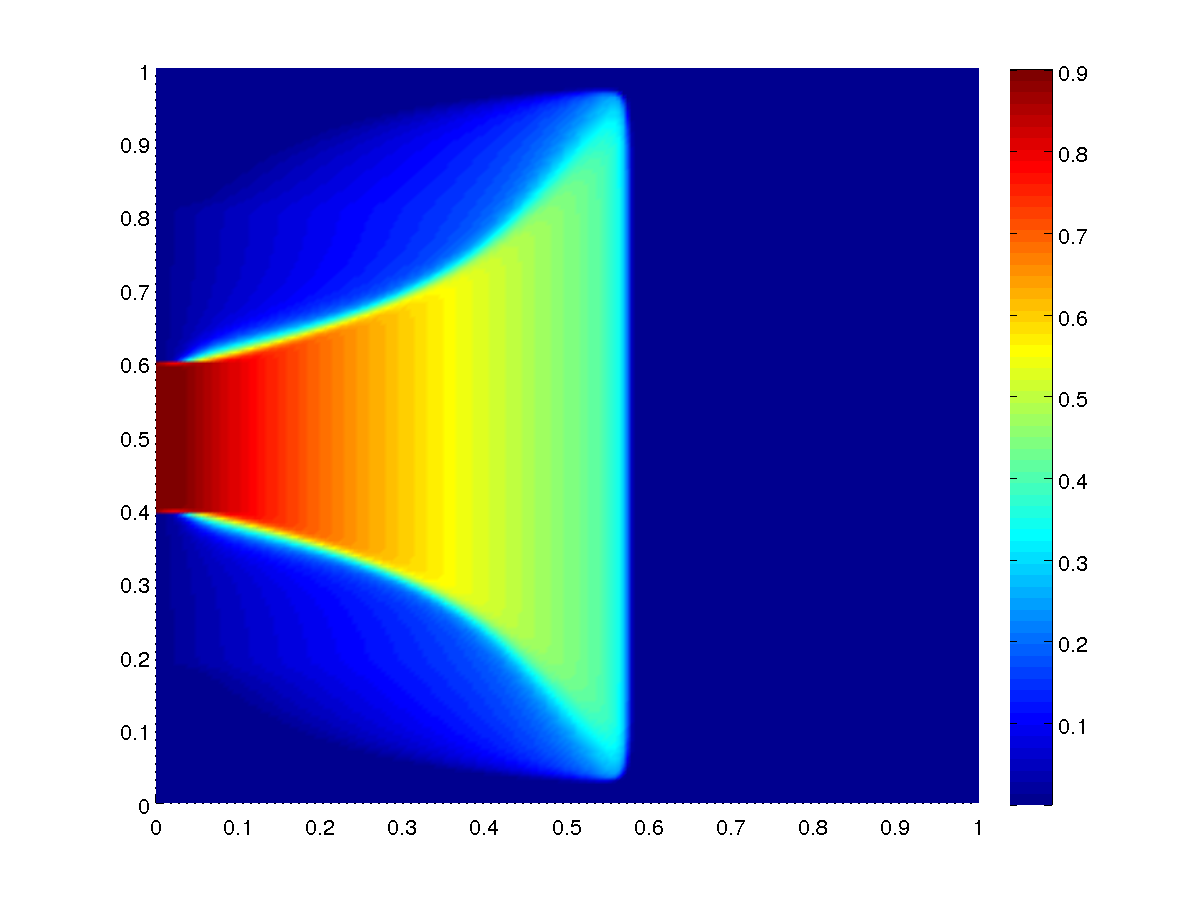}
\label{subfig:SM_ab_3_x_8}
}\hspace{-0.1cm}
\subfigure[finite-volume scheme]{
\includegraphics[scale=0.2]{./VEmodelS09-a5-200-200.png}
\label{subfig:SM_ab_3_NLT}
}
\caption{Numerical solution of the Multiscale model \eqref{eq:multiscale-dimensionless-a}, \eqref{eq:multiscale-dimensionless-b} using the multiscale algorithm \cite{Guo2014} in (a) and the finite-volume scheme in (b).}
\label{fig:VE-MS}
\end{figure}

\textbf{Example 1}: Figure \ref{fig:VE-MS} compares a numerical solution for the Multiscale model \eqref{eq:multiscale-dimensionless-a}, \eqref{eq:multiscale-dimensionless-b} (in the equivalent form of the VE-model) using the finite-volume scheme, presented in Section \ref{sec:finitevolumescheme}, with a numerical solution using the multiscale algorithm in \cite{Guo2014}. In both cases, a Cartesian grid with $200\times200$ elements is used, the end time is $T=0.3$ and the viscosity ratio is $M=5$. 

Figure \ref{fig:VE-MS} displays very similar discrete solutions. This indicates that both schemes have similar accuracy.

\textbf{Example 2}: Table \ref{table:time-2} compares the CPU-time for solving the Multiscale model using the finite-volume scheme, presented in Section \ref{sec:finitevolumescheme}, and using the algorithm in \cite{Guo2014}. In both cases, Cartesian grids with a fixed number of horizontal cells $N_x=100$ and a varying number of vertical cells $N_z=iN_x,\,i\in\{1,\,2,\,4,\,8,\,16\}$ are used. Table \ref{table:time-2} shows that the computational time using the finite-volume scheme is significantly reduced compared to that using the multiscale algorithm, mainly when the horizontal discretization is finer than in vertical one. This results from solving the elliptic equation of the vertically averaged pressure in the Multiscale model, which consumes most of the CPU-time. 

\begin{table}
\centering
\begin{tabular}{|c| c| c|}
\hline
Grid size $N_x\times N_z$ & Finite-volume scheme & Multiscale algorithm\\
\hline
$100 \times 100$ & 1.08 s & 0.91 s\\
$200 \times 100$ & 2.58 s & 2.97 s\\
$400 \times 100$ & 9 s & 12.32 s\\
$800 \times 100$ & 25 s & 64.85 s\\
$1600 \times 100$ & 109 s & 413 s\\ [0.5ex]
\hline
\end{tabular}
\caption{CPU-time for solving the VE-model using the finite-volume scheme compared to CPU-time for solving the Multiscale model using the multiscale algorithm with $T=0.3$ and viscosity ratio $M=5$.}
\label{table:time-2}
\end{table}

\textbf{Example 3}
Using the equivalence result from Theorem \ref{thm:VE-multi} and following Yortsos \cite{Yortsos}, we show in this example that the reconstruction function \eqref{eq:reconstruction} in domains with a small $\gamma$ represents strong gravity effects that lead to a fast segregation of the fluids, meaning that the vertical equilibrium assumption is fully satisfied. This, however, implies that the Multiscale model and the VI-model give results with similar accuracy. Consequently, the VI-model is sufficient to describe the fluids' flow as it has a lower computational complexity. 

We consider the Multiscale model \eqref{eq:multiscale1-a}, \eqref{eq:multiscale1-b} in a vertical cross-section $\Omega_v=(0,L)\times(0,H)$ of $\Omega\subset\mathbb{R}^3$
with negligible capillary forces and using the dimensionless variables \eqref{eq:dimensionlessvariables}, that is
\begin{align}
 -\partial_x \bigl( \hat{\lambda}_{tot}(\hat{S}_{i})\hat{\kappa}\partial_x\hat{p}\bigr)=0 \quad\text{ in }(0,1)\times(0,T),
 \label{eq:multiscale-dimensionless11-a}
\end{align}
and in $D\times(0,T)$
\begin{align}
\begin{array}{rl}
  p(x,z,t)=&\,\hat{p}(x,t)+ \pi_i(x,z,t),\\
  -\lambda_{tot}(S_i)\kappa \partial_x p=&u,\\
  \frac{1}{\gamma}\partial_z w=&-\partial_x u,\\
  f(S_i)w + \frac{g k (\rho_i-\rho_d)}{\gamma^2 q^3 \mu_d}\kappa f(S_i)\lambda_d(S_i)=&w_i ,\\
  \phi\partial_{t}S_i+\partial_x \big(uf(S_i)\big)+\frac{1}{\gamma}\partial_z\big( w_i \big)=&\, 0.\\
\end{array}
\label{eq:multiscale-dimensionless11-b}
\end{align}
The dimensionless reconstruction operator is now given as
\begin{align*}
 \pi_i(\cdot,z,\cdot)= -\dfrac{1}{\gamma q \mu_d}\int_0^z \Big(S_i(\cdot,r,\cdot)\rho_i + \big(1-S_i(\cdot,r,\cdot)\big)\rho_d\Big) g \,dr,
\end{align*}
for any $z\in (0,1)$.

To prove fluids' segregation, we consider the fourth equation in \eqref{eq:multiscale-dimensionless11-b}. With respect to the parameter $\gamma$, the term on the right side and the first term on the left side are of order $\mathcal{O}(1)$. However, the second term on the left side is of order $\mathcal{O}(\gamma^{-2})$, which diverges as $\gamma$ tends to zero. This leads to a contradiction unless the product $f\lambda_d$ satisfies
\begin{align*}
 f(S_i)\lambda_d(S_i)= 0,
\end{align*}
which implies either $\lambda_d(S_i)=0$ or $f(S_i)=0$, i.e. the fluids segregate. Another way to avoid the contradiction in the fourth equation in \eqref{eq:multiscale-dimensionless11-b} is to assume a small gravity effect by multiplying the gravity term in the first and the fourth equations in \eqref{eq:multiscale-dimensionless11-b} by a small parameter $N_G= \mathcal{O}(\gamma^2)$. By doing this, the second term in the fourth equation in \eqref{eq:multiscale-dimensionless11-b} is of order $\mathcal{O}(1)$. However, the gravity term in the first equation is of order $\mathcal{O}(\gamma)$. Therefore, this term can be neglected in media with a small parameter $\gamma$. 

\section{Models for Two-Phase Brinkman Flow in Flat Domains}
\label{sec:Brinkmanmodel}
In this part of the paper, we study two-phase flow in flat domains for Brinkman regimes. Examples of such regimes are media with high porosity. First, we present the main governing equations of the Brinkman two-phase flow model. Then, we derive a reduced model that describes two-phase Brinkman flow in flat domains. 

\subsection{The Brinkman Two-Phase Flow Model (BTP-Model)}
Brinkman's equations describe the flow of incompressible fluids through a porous medium with high porosity \cite{Brinkman1949}. Mathematically speaking, the Brinkman equations
\begin{equation}
\left(
\begin{array}{c}
-\mu_{e}\Delta u_{\alpha}+u_{\alpha}\\-\mu_{e}\Delta w_{\alpha}+w_{\alpha}
\end{array}\right)=-\dfrac{k_{r\alpha}}{\mu_{\alpha}}\textbf{K}(\nabla p_{\alpha}-\rho_{\alpha}\textbf{g}),
\label{eq:Darcy-perm1}
\end{equation}
are a regularization of Darcy's law. The terms $u_{\alpha}$ and $w_{\alpha}$ are the horizontal and vertical components of the velocity $\textbf{v}_{\alpha}$ for the phase $\alpha\in \{i,d\}$, respectively. The parameter $\mu_e>0$ is the effective viscosity. The second-order terms $\Delta u_{\alpha},\ \Delta w_{\alpha}$ on the left side of \eqref{eq:Darcy-perm1} correspond to components of the divergence of the stress tensor.

The validity of Brinkman's equation in porous media has been investigated by many researchers, see \cite{Auriault2009,Auriault2005} and the references therein. In media with high porosity, the effective viscosity $\mu_e$ is close to the dynamical viscosity. However, if the medium is macroscopiclly nonhomogeneous, the effective viscosity is very small ($\mu_e \ll 1$), such that the second order term in \eqref{eq:Darcy-perm1} is considered as a corrector of Darcy's law \cite{Auriault2005}. In addition, in porous media with almost unidirectional fluid flows, Brinkman's equations can be derived by averaging the Navier-Stokes equations over a representative averaging volume \cite{Hornung}. 

We define the vector $\textbf{V}_{\alpha}$, such that
\begin{align} 
\textbf{V}_{\alpha}=\left(\begin{array}{c}
U_{\alpha}\vspace{5pt}\\
W_{\alpha}
\end{array}\right)\coloneqq\left(
\begin{array}{c}
 -\mu_{e} \Delta u_{\alpha} +u_{\alpha}\vspace{5pt}\\
 -\mu_{e} \Delta w_{\alpha} +w_{\alpha}
\end{array}\right),
\label{eq:BrinkmanDarcy1}
 \end{align}
for each phase $\alpha\in \{i,d\}$. For incompressible fluids in a porous medium with the assumption of negligible gravity force, the Brinkman two-phase flow model (BTP-model) is a combination of the continuity equation, Brinkman's equations, and the incompressibility equation, that is
\begin{align}
\begin{array}{rl}
 \partial_{t} S_{\alpha} + \nabla\cdot \textbf{v}_{\alpha}&=0, \\
 \textbf{V}_{\alpha}&=-\dfrac{k_{r\alpha}}{\mu_{\alpha}} \textbf{K} \nabla p_\alpha,\\
   \nabla\cdot \textbf{v}&=0
\end{array} \quad\quad\text{ in }\Omega\times (0,T),
\label{eq:Brinkman-two-phase}
\end{align}
for both invading and defending fluids $\alpha\in\{i,d\}$, where $\Omega\coloneqq(0,L)\times(0,H)\subset \mathbb{R}^2$. Here, $\textbf{v}_{\alpha}=(u_{\alpha},w_{\alpha})$ is phase velocity, $\textbf{K}$ is permeability tensor, and $p_{\alpha}$ is phase pressure.
The third equation in \eqref{eq:Brinkman-two-phase} is the incompressibility relation for the total velocity $\textbf{v}=\textbf{v}_i+\textbf{v}_d$, which results from summing the continuity equation for the phases $\alpha\in\{i,d\}$ with the closure relation $S_{i}+S_{d}=1$. Applying the divergence operator $(\nabla\cdot)$ on $\textbf{V}$, then using the incompressibility equation, we obtain                                                                                                                                                                                                                                                                                                                                                                                                                                                                     
\begin{align}
 \nabla\cdot \textbf{V}  = \nabla\cdot\left(\begin{array}{c}
U\vspace{5pt}\\ W\end{array}\right)& =  -\mu_e \nabla\cdot\left(\begin{array}{c}\Delta u\\ \Delta w \end{array}\right)+ \nabla\cdot\left(\begin{array}{c} u\\  w \end{array}\right),\nonumber \\ &= -\mu_e\,\Delta (\nabla\cdot \textbf{v}) +\nabla \cdot \textbf{v}=0.
\label{eq:incompressibilty2}
\end{align}
Model \eqref{eq:Brinkman-two-phase} is completed with the initial and boundary conditions in \eqref{eq:TPF-IBC}, with an extra periodic condition on the Laplacian $\Delta w$ at the impermeable boundary $\partial_{\text{imp}}\Omega$,
\begin{align}
\begin{array}{rll}
 S_i(\cdot,\cdot,0)&=S_{0}, &\text{ in } \Omega, \\
 S_i&=S_{\text{inflow}}, &\text{ on } \partial\Omega_{\text{inflow}}\times [0,T],\\
 \textbf{n}\cdot\textbf{v}_{\alpha}&= 0,  &\text{ on }\partial \Omega_{\text{imp}}\times [0,T],\\
 \Delta w_{\alpha}(x,0,t)&=\Delta w_{\alpha}(x,1,t),   &\text{ for }x\in[0,L],\,t\in [0,T].
\end{array}
\label{eq:BVE-IBC}
\end{align}
Here, $\textbf{n}$ is the outer normal vector of the boundary $\partial \Omega_{\text{imp}}$. Note that the extra periodic condition in \eqref{eq:BVE-IBC} has been added to account for the new higher-order terms in the model.
\subsection{The Brinkman Vertical Equilibrium Model}
\label{sec:brinkman}     
In this section, we derive the Brinkman Vertical Equilibrium Model (BVE-model) following the analysis in \cite{Yortsos}. First, we rescale the variables in the BTP-model \eqref{eq:Brinkman-two-phase} into dimensionless ones, which leads to a dimensionless model of two unknowns depending on the geometrical parameter $\gamma=\frac{H}{L}$. Second, we apply formal asymptotic analysis to the dimensionless model for $\gamma$ tends to $0$. This produces a nonlocal evolution equation for the remaining unknown saturation. The higher-order Brinkman terms lead to  a third-order term of mixed derivatives.

\subsubsection{The Dimensionless BTP-Model}
The first step to derive the BVE-model is rescaling the BTP-model \eqref{eq:Brinkman-two-phase} using the dimensionless variables \eqref{eq:dimensionlessvariables}. Applying the chain rule to \eqref{eq:BrinkmanDarcy1}, then defining the parameters
\begin{equation}
 \beta_x \coloneqq \dfrac{\mu_e}{L^2},\quad \beta_z\coloneqq\dfrac{\mu_e}{H^2},
 \label{eq:betax-y}
\end{equation}
and the dimensionless components
\begin{align}
 \overline{U}_{\alpha}\coloneqq\dfrac{U_{\alpha}}{q}, \quad\quad \overline{W}_{\alpha}\coloneqq\dfrac{W_{\alpha}}{q},
 \label{eq:B-dime-compo}
\end{align}
yield
\begin{align}
\begin{array}{ll}
 \overline{U}_{\alpha}=&\overline{u}_{\alpha}-\beta_x\partial_{\overline{x}\overline{x}}\overline{u}_{\alpha}  -\beta_z\partial_{\overline{z}\overline{z}}\overline{u}_{\alpha},\vspace{0.1cm}\\
 \overline{W}_{\alpha}=&\overline{w}_{\alpha}-\beta_x\partial_{\overline{x}\overline{x}}\overline{w}_{\alpha}  -\beta_z\partial_{\overline{z}\overline{z}}\overline{w}_{\alpha}.
\end{array}
 \label{eq:operatorNon}
\end{align}
Applying the chain rule to equations \eqref{eq:Brinkman-two-phase} and \eqref{eq:incompressibilty2}, using equations \eqref{eq:dimensionlessvariables} and \eqref{eq:B-dime-compo}, then omitting the bar-sign leads to
\begin{align}
\begin{array}{rl}
 \partial_{t}S_{\alpha} +\partial_{x}u_{\alpha}+(1/\gamma) \partial_{z}w_{\alpha}&=0,\vspace{3pt}\\
 U_{\alpha}&=-\lambda_{\alpha}(S_{\alpha})\, \kappa_x\, \partial_{x} p_{\alpha},\vspace{3pt}\\
 (\gamma/\sigma) W_{\alpha}&=-\lambda_{\alpha}(S_{\alpha})\, \kappa_z\, \partial_{ z}  p_{\alpha},\vspace{3pt}\\
  \partial_{x}u+(1/\gamma) \partial_{z}w &=0,\vspace{3pt}\\
   \partial_{x}U+(1/\gamma) \partial_{z}W &=0
\end{array}
\label{eq:B-Dimensionlesstwo-phase2}
\end{align}
in $D\times(0,T)$, for both invading and defending phases $\alpha\in\{i,d\}$, where $D\coloneqq(0,1)\times(0,1)$ is the dimensionless spatial domain. Here $\lambda_{\alpha}=\tfrac{\mu_d k_{r\alpha}}{\mu_{\alpha}}$ is the dimensionless mobility, 
 \begin{align}
  \gamma\coloneqq\dfrac{H}{L} \quad \text{and} \quad \sigma\coloneqq\dfrac{k_z}{k_x}.
  \label{eq:gamma-segma}
 \end{align}
The assumption of negligible capillary pressure implies $p_i=p_d=:p$ and the phases' velocities satisfy 
\begin{align}
  U_{\alpha}=f(S_{\alpha})U,\quad\quad  W_{\alpha}=f(S_{\alpha})W.
  \label{eq:Dimensionlessfrac}
\end{align}
The dimensionless model \eqref{eq:B-Dimensionlesstwo-phase2} in $D\times(0,T)$ is now summarized such that the unknown variables $S_i,\,p,\,U$, and $W$ are associated with the geometrical parameter $\gamma$,
\begin{align}
\begin{array}{rl}
\partial_{t}S_i^{\gamma}-\beta_x\partial_{xxt}S_i^{\gamma}-\beta_z \partial_{zzt}S_i^{\gamma} &+\partial_{x}\left(f(S_i^{\gamma})U^{\gamma}\right)\\&+ \dfrac{1}{\gamma}\partial_{z}\left(f(S_i^{\gamma})W^{\gamma}\right)=0,\\
U^{\gamma}&=-\lambda_{tot}(S_i^{\gamma}) \kappa_x \partial_{x} p^{\gamma},\\
 \gamma  W^{\gamma}&=-\lambda_{tot}(S_i^{\gamma}) \kappa_z \partial_{z} p^{\gamma},\\
  \partial_{x}U^{\gamma}+\dfrac{1}{\gamma} \partial_{z}W^{\gamma} &=0.
\end{array}
\label{eq:B-Dimensionlessfractionalformulation}
\end{align}
\subsubsection{Asymptotic Analysis}
The second step in deriving the BVE-model is applying formal asymptotic analysis, with respect to $\gamma$, to the dimensionless BTP-model \eqref{eq:B-Dimensionlessfractionalformulation}. The existence of a weak solution $(S_i^{\gamma},p^{\gamma},U^{\gamma}, W^{\gamma})$ for model \eqref{eq:B-Dimensionlessfractionalformulation} with appropriate initial and boundary conditions is established in \cite{Coclite2014}. We assume in this section that each component in $(S_i^{\gamma},p^{\gamma},U^{\gamma}, W^{\gamma})$ is smooth and can be written in terms of the asymptotic expansions
\begin{align}
 \begin{array}{ll}
  Z^{\gamma}=Z_{0}+\gamma Z_{1}+\mathcal{O}(\gamma^2),\,& Z^{\gamma}\in\{S_i^{\gamma},\,p^{\gamma},\,U^{\gamma}\},\vspace{5pt}\\
  W^{\gamma}=\gamma W_{1}+\mathcal{O}(\gamma^2).
 \end{array}
 \label{eq:B-asymptoticexpansion}
\end{align}
Note that the second equation in \eqref{eq:B-asymptoticexpansion} corresponds to the assumption that the vertical velocity in flat domains is small.

Using the asymptotic expansion of $S^{\gamma}$ in \eqref{eq:B-asymptoticexpansion} and Assumption \ref{ass:assumptions}, we have the Taylor expansions
\begin{align}
\begin{array}{cl}
G(S_i^{\gamma})&=G(S_{i,0})+ G'(S_{i,0})(\gamma S_{i,1})+\mathcal{O}(\gamma^{2}),
\end{array}
 \label{eq:B-asymptoticexpansion1}
\end{align}
for $G\in\{\lambda_{tot},\,f\}$. The first incompressibility relation in \eqref{eq:B-Dimensionlessfractionalformulation} allows to write the continuity equation in nonconservative form. Moreover, Assumption \ref{ass:assumptions}(1) implies that the parameter $\sigma$, defined in \eqref{eq:gamma-segma}, satisfies $\sigma=1$ and, therefore, we set $\kappa=\kappa_x=\kappa_z$. Substituting equation \eqref{eq:B-asymptoticexpansion} and \eqref{eq:B-asymptoticexpansion1} into \eqref{eq:B-Dimensionlessfractionalformulation} and neglecting the capillary pressure terms in the third and fourth equations of \eqref{eq:B-Dimensionlessfractionalformulation}, the terms of order $\mathcal{O}(1)$ satisfy
\begin{align}
\begin{array}{rl}
 \partial_{t}S_{i,0}-\beta_x\partial_{xxt}S_{i,0}&-~\beta_z \partial_{zzt}S_{i,0} +\partial_{x}\big(f(S_{i,0})U_{0}\big)\vspace{5pt}\\&+ \partial_{z}\big(f(S_{i,0})W_1\big)=\mathcal{O}(\gamma),\vspace{5pt}\\
  U_{0}&=-\lambda_{tot}(S_{i,0}) \kappa \partial_{x} p_{0},\vspace{5pt}\\
  \lambda_{tot}(S_{i,0}) \kappa \partial_{z} p_{0}&=\mathcal{O}(\gamma),\vspace{5pt}\\
 \partial_{x}U_{0}+\partial_{z}W_{1} &=\mathcal{O}(\gamma).
\end{array}
\label{eq:B-asymptotic}
\end{align}
Using the positivity of the permeability $\kappa$ and the total mobility $\lambda_{tot}$ (see Assumptions \ref{ass:assumptions}(1) and \ref{ass:assumptions}(3)), the third equation of \eqref{eq:B-asymptotic} implies that $p_0$ is independent of the $z$-coordinate,
\begin{equation}
 p_{0}=p_{0}(x,t).
 \label{eq:B-pressure-y-indep}
\end{equation}
This result on $p_0$ substitutes the vertical equilibrium assumption \eqref{eq:ref-pressure}, which is however a consequence of the expansions \eqref{eq:B-asymptoticexpansion} and the scaling in \eqref{eq:B-Dimensionlessfractionalformulation}. 

Integrating the last equation in \eqref{eq:B-asymptotic} over the vertical direction from $0$ to $1$ then using the zero-Neumann boundary condition and the periodicity of $\Delta w_{\alpha}$ at the impermeable boundary $\partial_{\text{imp}} D$ in equation \eqref{eq:BVE-IBC}, we obtain
\begin{eqnarray*}
\partial_{x}\int_{0}^{1}U_{0}\,dz =-\int_{0}^{1}\partial_{z}W_{1}\,dz=-\int_{0}^{1}\partial_z (w_1-\mu_e \Delta w_1) \,dz=0.
\end{eqnarray*}
This yields 
\begin{align}
  \int_{0}^{1}U_{0}\,dz= q(t),
  \label{eq:U1}
\end{align}
for some positive function $q=q(t)$ with $q(t)>0$ for all $t\in[0,T]$. Substituting the second equation in \eqref{eq:B-asymptotic} into equation \eqref{eq:U1} yields
\begin{equation*}
-\int_{0}^{1}\lambda_{tot}(S_{i,0}) \kappa \partial_{x} p_{i,0}\,dz=q(t).
\end{equation*}
Then, using equation \eqref{eq:B-pressure-y-indep}, we have
\begin{equation}
\partial_{x} p_{0}(x,t)=-\dfrac{q(t)}{\int_{0}^{1}\lambda_{tot}(S_{i,0}(x,z,t)) \kappa(x,z)\,dz},
\label{eq:U2}
\end{equation}
for all $x\in (0,1)$ and $t\in(0,T)$. Substituting \eqref{eq:U2} into the second equation in \eqref{eq:B-asymptotic}, we obtain a nonlocal saturation-dependent formula for $U_0$,
\begin{equation}
 U_{0}[x,z,t;S_{i,0}]=\dfrac{q(t)\lambda_{tot}\bigl(S_{i,0}\bigr)\kappa(x,z) }{\int_{0}^{1}\lambda_{tot}\bigl(S_{i,0}\bigr) \kappa(x,z)\,dz},
  \label{eq:U3}
\end{equation}
for all $(x,z)\in D$ and $t\in(0,T)$. Consequently, the incompressibility relation in \eqref{eq:B-asymptotic} yields also a nonlocal saturation-dependent formula for $W_1$,
\begin{equation}
 W_1[x,z,t;S_{i,0}]=-\partial_{x}\int_{0}^{z}U_0[x,r,t;S_{i,0}]\,dr,
 \label{eq:U4}
\end{equation}
for all $(x,z)\in D$ and $t\in(0,T)$.
Using equation \eqref{eq:U3} and \eqref{eq:U4}, omitting the subscripts $\{0,1\}$, rescaling the time using 
\begin{align}
 t\mapsto \bar{t} = \int_0^t q(r)\,dr+q(0)t,
 \label{eq:B-timerescale}
\end{align}
then omitting the subscripts, system \eqref{eq:B-asymptotic} reduces to a third-order pseudo-parabolic equation of saturation alone that accounts to the vertical dynamics of the fluids. We call this equation the Brinkman Vertical Equilibrium model (BVE-model),
\begin{align}
 \partial_{t}S_i+\partial_{x}\left(f(S_i)U[\cdot,\cdot;S_i]\right)+&\partial_{z}\left(f(S_i)W[\cdot,\cdot;S_i]\right)\nonumber\\&-\beta_x\partial_{xxt}S_i-\beta_z \partial_{zzt}S_i  =0,
\label{eq:Brinkman1}
\end{align}
in $D\times (0,T)$ where we have for all $z\in(0,1)$
\begin{align}
\begin{array}{rl}
 U[\cdot,\cdot;S_i]&=\dfrac{\lambda_{tot}(S_i)\kappa}{\int_{0}^{1}\lambda_{tot}(S_i) \kappa \,dz}, \vspace{5pt}\\
 W[\cdot,z;S_i]&=-\partial_{x}\int_{0}^{z}U(\cdot,r;S_i(\cdot,r,\cdot))\,dr .
 \end{array}
 \label{eq:Brinkman2}
\end{align}

\begin{remark}
 In asymptotically flat domains, the definition of the parameters $\beta_x$ and $\beta_z$ yields 
\begin{align*}
\begin{array}{ll}
 \beta_x=\mathcal{O}(\gamma^2),\quad\quad &\beta_z=\mathcal{O}(1).
\end{array}
\end{align*} 
This implies that the dissipative effect due to viscosity in the horizontal direction is negligible with respect to that in the vertical direction. Then, the BVE-model \eqref{eq:Brinkman1}, \eqref{eq:Brinkman2} simplifies to
\begin{align*}
 \partial_{t}S_i+\partial_{x}\left(f(S_i)U[\cdot,\cdot;S_i]\right)+&\partial_{z}\left(f(S_i)W[\cdot,\cdot;S_i]\right)-\beta_z \partial_{zzt}S_i  =0,
\end{align*}
where $U$ and $W$ satisfy \eqref{eq:Brinkman2}.
\end{remark}

\begin{remark}
 Except of the nonlocal definition of the velocity components in \eqref{eq:Brinkman2}, it is notable that the BVE-model \eqref{eq:Brinkman1} has a similar structure as the pseudo-parabolic model from \cite{HG1993} that has been mathematically analyzed in \cite{Cao2015,FanPop2011}.
\end{remark}

\section{Comparison of the Brinkman-Type Models in Flat Domains}
\label{sec:B-numericalresults}

This section, we verify the limit $\gamma\rightarrow 0$ in the BTP-model \eqref{eq:Brinkman-two-phase} numerically. We also demonstrate the computational efficiency of the BVE-model \eqref{eq:Brinkman1}, \eqref{eq:Brinkman2} compared to the BTP-model \eqref{eq:B-Dimensionlessfractionalformulation}. In addition, some qualitative properties of the solutions of the BVE-model are shown by comparing it with the VE-model \eqref{eq:VEmodel}. To simplify the notation, we set here $S \coloneqq S_i$. For the numerical examples, we consider the BVE-model with a small capillary force such that
\begin{align}
\partial_{t}S & +\partial_{x}\left(f(S)U[\cdot,\cdot;S]\right)+\partial_{z}\left(f(S)W[\cdot,\cdot;S]\right)-\epsilon_x\partial_x\left(H(S)\partial_x S\right)\nonumber\\&- \epsilon_z\partial_z\left(H(S)\partial_z S\right)-\beta_x \partial_{xxt}S- \beta_z \partial_{zzt}S=0,
\label{eq:Bmodel1}
\end{align}
with
\begin{equation}
\begin{array}{rl}
U[\cdot,\cdot;S]&=\dfrac{\lambda_{tot}(S)\kappa}{\int_{0}^{1}\lambda_{tot}(S) \kappa\, dz},\vspace{5pt}\\
W[\cdot,z;S]&= -\partial_{x}\int_{0}^{z}U[\cdot,r;S(\cdot,r,\cdot)]\,dr,
\end{array}
\label{eq:Bmodel2}
\end{equation} 
for all $z\in(0,1)$. We set $\epsilon_x=\sqrt{\beta_x}$, $\epsilon_z=\sqrt{\beta_z}$ and choose the nonlinear diffusion function $H=H(S)$ such that
\begin{align}
 H(S)=\kappa \frac{MS^2(1-S)^2}{MS^2+(1-S)^2},
 \label{eq:Bmodel3}
\end{align}
where $M$ is the viscosity ratio of the defending phase and the invading phase.

Dividing the domain $D=(0,1)\times(0,1)$ into $N_z\in\mathbb{N}$ layers, each of width $\Delta z=\frac{1}{N_z}$, the BVE-model \eqref{eq:Bmodel1}, \eqref{eq:Bmodel2} reduces into a coupled system of $N_z$ equations. Applying Euler's method to the time derivative, a finite volume discretization for \eqref{eq:Bmodel1} is given by
\begin{align}
\hspace{-0.2cm}
\begin{array}{ll}
& \frac{S_{i,j}^{n+1}-S_{i,j}^{n}}{\Delta t}\hspace{-0.1cm} = -\dfrac{1}{\Delta x\Delta z}\sum_{l\in\theta_{i,j}}\mathcal{F}_{l}\left(\bar{S}_{i,j}^{n},S_{i,j}^{n},S_{(i,j)_l}^{n}\right)\vspace{3pt}\\
&+  \frac{\epsilon_x}{\Delta x^2}\left( H(S_{i+\frac{1}{2},j}^n)( S_{i+1,j}^n -S_{i,j}^n)-H(S_{i-\frac{1}{2},j}^n)( S_{i,j}^n -S_{i-1,j}^n)\right)\vspace{3pt} \\
&+ \frac{\epsilon_z}{\Delta z^2}\left( H(S_{i,j+\frac{1}{2}}^n)( S_{i,j+1}^n -S_{i,j}^n)-H(S_{i,j-\frac{1}{2}}^n)( S_{i,j}^n -S_{i,j-1}^n)\right)\vspace{3pt}\\
& + \frac{\beta_x }{\Delta x^2\Delta t}\Big( ( S_{i+1,j}^{n+1}-2S_{i,j}^{n+1}+S_{i-1,j}^{n+1} )-(S_{i+1,j}^n-2S_{i,j}^n+S_{i-1,j}^n)\Big)\vspace{3pt} \\
& + \frac{\beta_z }{\Delta z^2\Delta t}\Big(( S_{i,j+1}^{n+1}-2S_{i,j}^{n+1}+S_{i,j-1}^{n+1} )-(S_{i,j+1}^n-2S_{i,j}^n+S_{i,j-1}^n )\Big) .
\end{array}
\label{eq:B-numericalscheme1}
\end{align}
The numerical flux function and the discrete velocity vector $\textbf{V}_{l}^{n}=(U_{l}^{n},W_{l}^{n})^{T}$ are defined as in Section \ref{sec:finitevolumescheme}. The initial condition $S_{i,j}^0$ is defined in each cell $T_{i,j}$ as follows
\begin{align*}
 S^0_{i,j}=\dfrac{ (1-i\Delta x)^2S_{\text{inflow},j}}{10^5 (i\Delta x)^2+(1-i\Delta x)^2},
\end{align*}
for a constant $S_{\text{inflow},j} \in[0,1]$, $j=1,\cdots,N_z$.

The numerical scheme \eqref{eq:B-numericalscheme1} satisfies the following two stability criteria:

  	\textbf{Mass conservation}:
	 Thanks to the definition of the numerical flux in \eqref{eq:numericalfluxfunction} and the symmetry of the discrete Laplacian operator, the scheme \eqref{eq:B-numericalscheme1} is mass-conservative, i.e., if $E_l=E_{l'}$ such that $l\in \theta_{i,j}$ and $l'\in \theta_{(i,j)_l}$, then, for all $P^n_{i,j}\in[0,1]$, it holds that
	\begin{equation*}
	  F_l(\bar{P}_{i,j}^n, P_{i,j}^n,P_{(i,j)_l}^n)=-F_{l'}(\bar{P}^n_{(i,j)_l},P^n_{(i,j)_l},P^n_{i,j}),
	\end{equation*}
	where 
	\begin{align*}
	 &F_l(\bar{P}^n_{i,j}, P^n_{i,j},P^n_{(i,j)_l})= -\dfrac{1}{\Delta x\Delta z}\mathcal{F}_{l}\left(\bar{S}_{i,j}^{n},S_{i,j}^{n},S_{(i,j)_l}^{n}\right)\\&+ \frac{\epsilon_l | E_l|^2}{\Delta x^2\Delta z^2} H(P_l)\big(P_{(i,j)_l} - P_{i,j} \big)-\frac{\beta_l | E_l|^2}{\Delta x^2 \Delta z^2 \Delta t}\big( P_{(i,j)_l} - P_{i,j}\big).
	\end{align*}

	\textbf{Discrete incompressibility}:
	If the vector of discrete velocity \emph{$\textbf{v}_{l}^{n}[\bar{S}_{i,j}^{n}]=(u_{l}^{n}[\bar{S}_{i,j}^{n}],~w_{l}^{n}[\bar{S}_{i,j}^{n}])^{T}$}, $l\in\theta_{i,j}$, is defined as in \eqref{eq:discretehorizontalvelocity1} and \eqref{eq:discretehorizontalvelocity2}, then \emph{$\sum_{l\in\theta_{i,j}}\textbf{n}_{l}\cdot \textbf{v}_{l}^{n}[\bar{S}_{i,j}^{n}]=0$}.

\begin{figure}
\centering
\subfigure[BTP-model for $\gamma=1$]{
\includegraphics[scale=0.25]{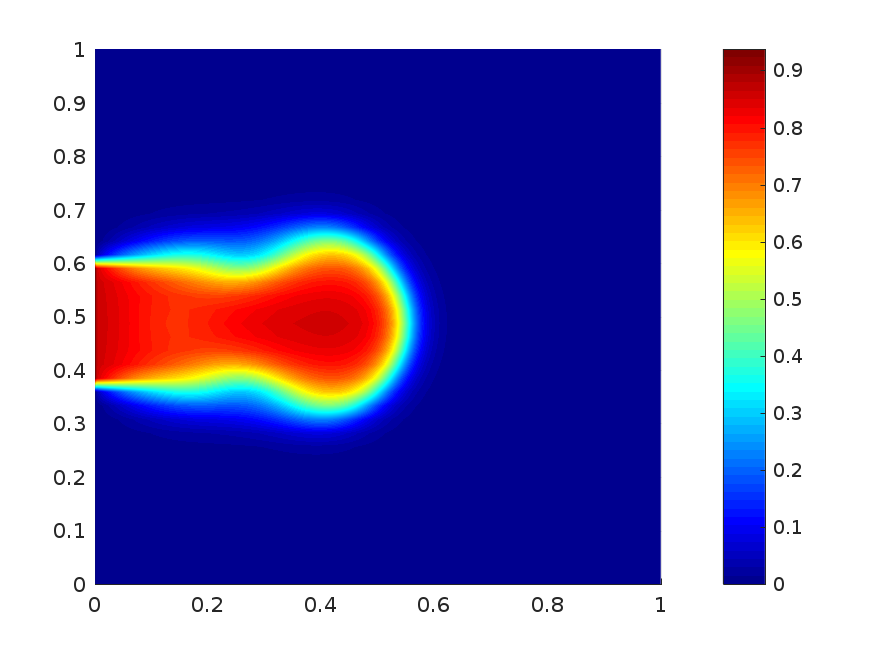}
\label{subfig:BTP-L5}
}\hspace{-0.1cm}
\subfigure[BTP-model for $\gamma=1/4$]{
\includegraphics[scale=0.25]{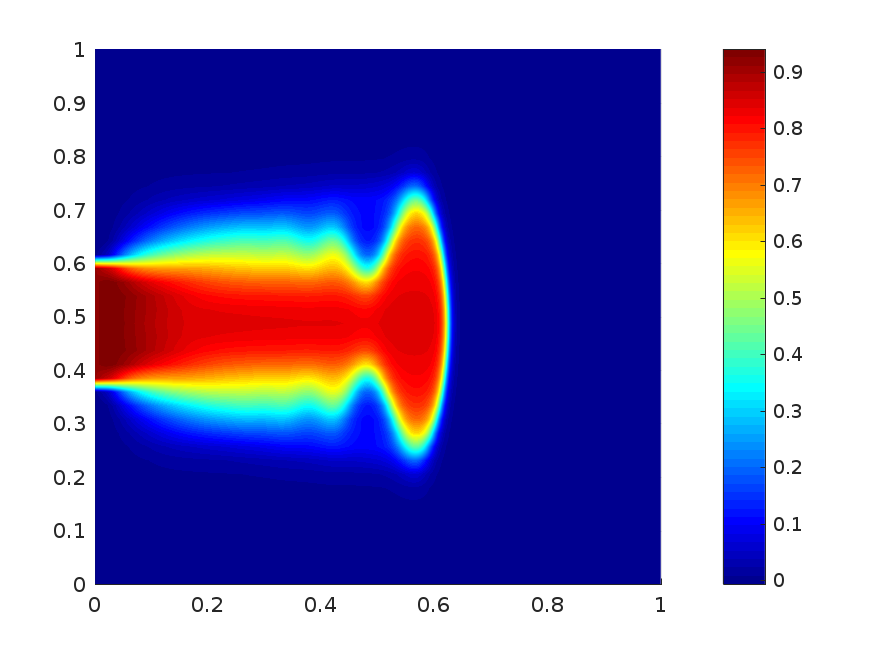}
\label{subfig:BTP-L20}
}\\
\subfigure[BTP-model for $\gamma=1/16$]{
\includegraphics[scale=0.25]{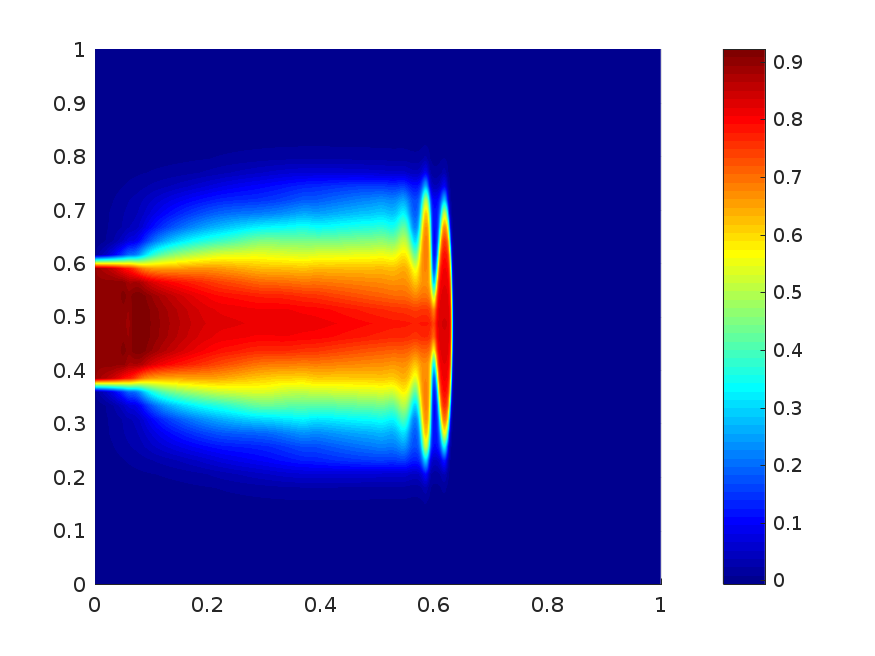}
\label{subfig:BTP-L80}
}\hspace{-0.1cm}
\subfigure[BTP-model for $\gamma=1/64$]{
\includegraphics[scale=0.25]{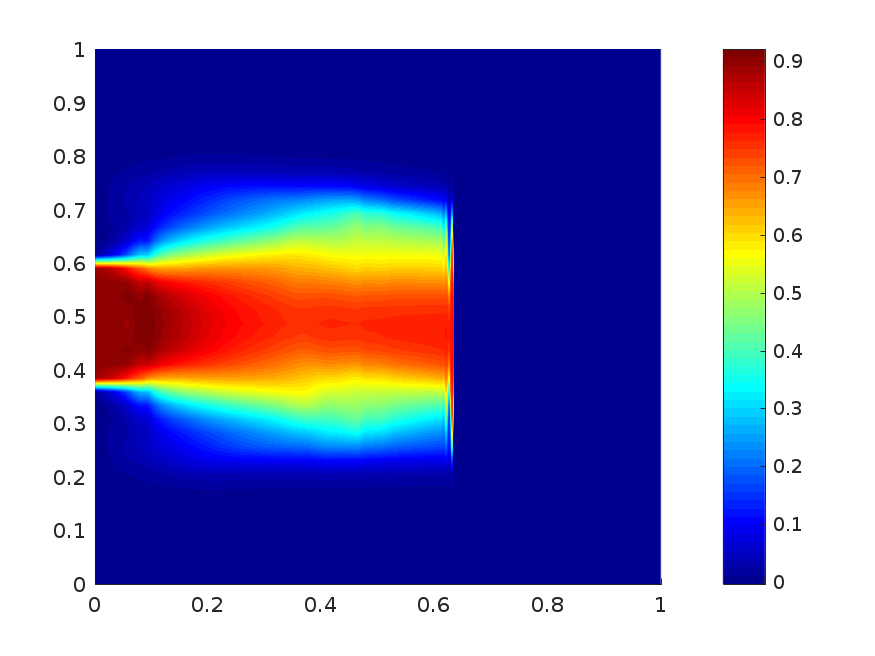}
\label{subfig:BTP-L480}
}\\
\subfigure[BVE-model]{
\includegraphics[scale=0.25]{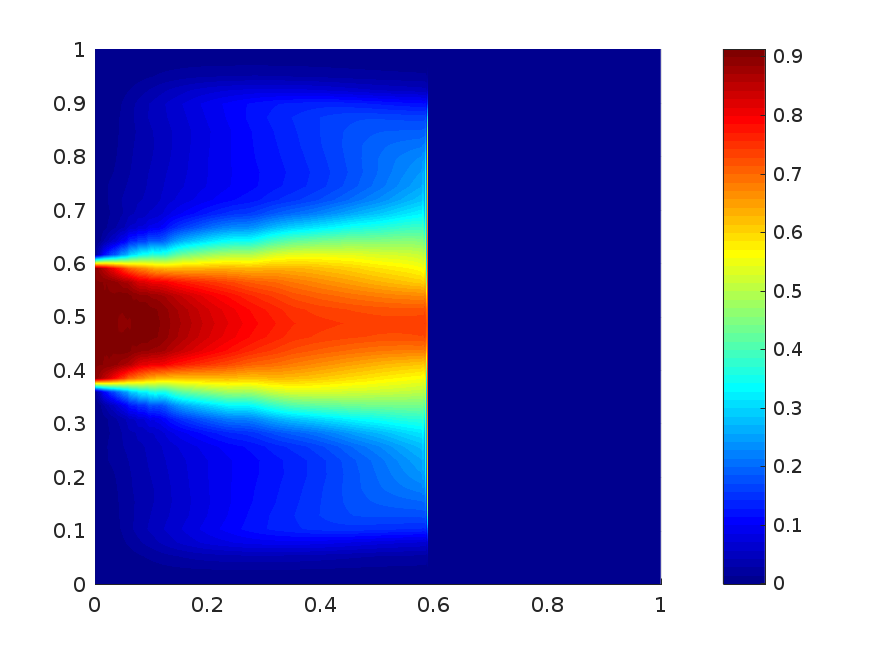}
\label{subfig:BVE}
}\hspace{-0.2cm}
\subfigure[Middle layer]{
\includegraphics[scale=0.28]{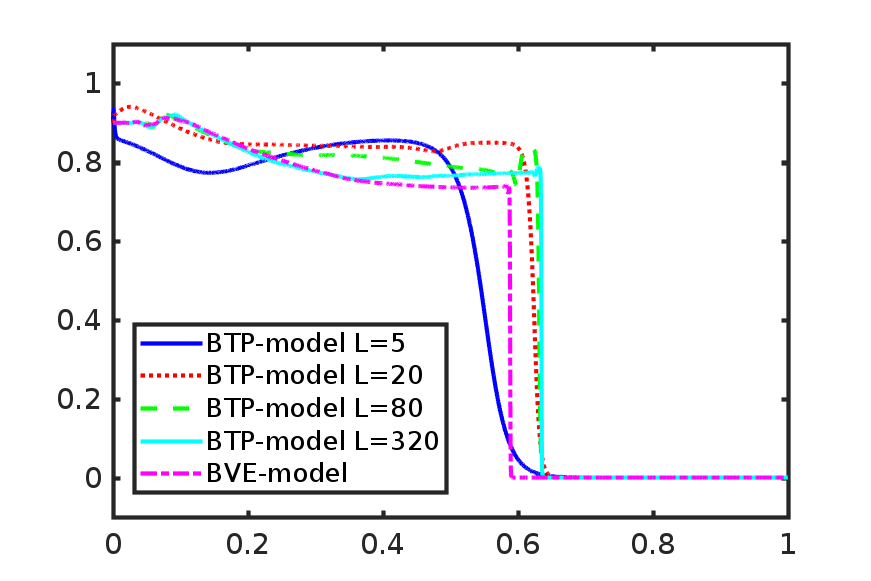}
\label{subfig:1DBcomparison1}
}
\caption{Numerical solutions of the BTP-model \eqref{eq:B-Dimensionlessfractionalformulation22} with $\gamma\in\{1,\,\frac{1}{4},\,\frac{1}{16},\,\frac{1}{64}\}$ in (a)-(d), the BVE-model \eqref{eq:Bmodel1}, \eqref{eq:Bmodel2} in (e) and the saturation profiles in the middle layer in (f), for $\mu_e=10^{-2}$ and $T=0.5$.}
\label{fig:Bcomparison1}
\end{figure}
	
\subsection{The BVE-Model vs. the BTP-Model}
In this section, we verify the limit $\gamma\rightarrow 0$ in the BTP-model \eqref{eq:Brinkman-two-phase} numerically. Furthermore, we investigate the computational efficiency of the BVE-model \eqref{eq:Bmodel1}, \eqref{eq:Bmodel2} compared to the BTP-model \eqref{eq:Brinkman-two-phase}. We consider the BTP-model \eqref{eq:Brinkman-two-phase} with the nonlinear diffusion function $H(S)$ such that
\begin{align}
\partial_{t}S^{\gamma}-\beta_x\partial_{xxt}S^{\gamma}-\beta_z\partial_{zzt}S^{\gamma} +&\partial_{x}\left(f(S^{\gamma})U^{\gamma}\right)+ \dfrac{1}{\gamma}\partial_{z}\left(f(S^{\gamma})W^{\gamma}\right)\nonumber\\+\epsilon_x\partial_x\big(H(S^{\gamma})\partial_x S^{\gamma}\big) & + ~ \epsilon_z \partial_{z}\big(H(S^{\gamma})\partial_{z}S^{\gamma}\big) =0,\nonumber\\
U^{\gamma}&=-\lambda_{tot}(S^{\gamma}) \kappa \partial_{x} p^{\gamma},\nonumber\\
 \gamma  W^{\gamma}&=-\lambda_{tot}(S^{\gamma}) \kappa \partial_{z} p^{\gamma},\nonumber\\
  \partial_{x}U^{\gamma}+\dfrac{1}{\gamma} \partial_{z}W^{\gamma} &=0
\label{eq:B-Dimensionlessfractionalformulation22}
\end{align}
in $D\times (0,T)$, where $\gamma=\tfrac{H}{L}$ is the geometrical parameter of the domain. For solving this model, we use the IMPES-method \cite{HuberHelmig1999} with pressure $p=1$ at the inflow boundary and $p=0$ at the outflow boundary.

In the following we present two examples: the first example shows the convergence of the numerical solution for the BTP-model \eqref{eq:B-Dimensionlessfractionalformulation22} to the corresponding solution of the BVE-model as the geometrical parameter $\gamma$ tends to $0$. The second example demonstrates the computational efficiency of the BVE-model over that of the BTP-model.

\textbf{Example 1}: Consider the viscosity ratio $M=2$ and the inflow condition
\begin{align}
 S_{\text{inflow}}(z)=\left\{ \begin{array}{cll}
                           0  \quad &: & z\leq \frac{2}{5} \text{ and } z>\frac{3}{5},\\
			  0.9 \quad &: &\frac{2}{5}<z\leq \frac{3}{5}.
                          \end{array}\right.
 \label{eq:inflowcond}
\end{align}
Figures (\ref{subfig:BTP-L5} - \ref{subfig:BTP-L480}) present the numerical solutions for the BTP-model \eqref{eq:B-Dimensionlessfractionalformulation22}, using different geometrical parameters $\gamma\in \{1, \frac{1}{4},\tfrac{1}{16},\,\tfrac{1}{64}\}$, respectively, that correspond to domains with a fixed width $H=5$ and a varying length $L\in\{5,\,20,\,80,\,320\}$. Figure \ref{subfig:BVE} presents the numerical solution of the BVE-model \eqref{eq:Bmodel1}, \eqref{eq:Bmodel2}. The parameters $\beta_x$ and $\beta_z$ (see \eqref{eq:betax-y}) in the BVE-model are chosen such that $\mu_e=10^{-2}$, $H=5$ and $L=320$. The spatial domain $D=(0,1)\times(0,1)$ is discretized for both models into a uniform Cartesian grid with $N_x=8000$ cells in the horizontal direction and $N_z=40$ cells in the vertical and the end time is $T=0.5$. 

The results in Figure \ref{fig:Bcomparison1} suggest that numerical solutions of the BTP-model \eqref{eq:B-Dimensionlessfractionalformulation22} converge to the numerical solution of the BVE-model as the geometrical parameter $\gamma$ tends to $0$. Moreover, the spreading speed of the invading fluid using the BVE-model is captured very well. For a better presentation of the convergence, Figure \ref{subfig:1DBcomparison1} presents the numerical solutions for both models in the middle layer of the domain. It shows that the saturation distribution and the spreading speed of the BTP-model with $\gamma=1/64$ match with those of the BVE-model.
 
Example 1 represents a validation criterion on the BVE-model due to two related reasons: firstly, the convergence of solutions for the BTP-model to that of the BVE-model as $\gamma$ tends to zero. Secondly, the spreading speed using the BVE-model matches with that of the BTP-model in flat domains.

\begin{table}
\centering
\begin{tabular}{|c | c | c | c | c}
\hline
$N_z$ & $N_x$ & Model & CPU-time  \\ [0.3ex]
\hline
 & 100 & BVE-model &  0.75 s\\
10 & 100 & BTP-model &  1.08 s\\
 & 1000 & BVE-model &  26.7 s\\
 & 1000 & BTP-model &  39.4 s\\
\hline
 & 100 & BVE-model &  8.2 s\\
80 & 100 & BTP-model &  5.53 s\\
 & 1000 & BVE-model & 462.26 s\\
 & 1000 & BTP-model &  678 s\\
\hline
\end{tabular}
\label{table:nonlin}
\caption{Comparison of CPU-time for the BVE-model \eqref{eq:Bmodel1}, \eqref{eq:Bmodel2} with the BTP-model \eqref{eq:B-Dimensionlessfractionalformulation22} using $\gamma=100$ and $\beta_x=\beta_z=10^{-6}$.}
\label{table:BVE-BTP}
\end{table}

\textbf{Example 2}: Consider the inflow condition in \eqref{eq:inflowcond} and set $\gamma=100$. In Table \ref{table:BVE-BTP}, we compare the CPU-time required to solve the BVE-model \eqref{eq:Bmodel1}, \eqref{eq:Bmodel2} to that required to solve the BTP-model \eqref{eq:B-Dimensionlessfractionalformulation22}. From Table \ref{table:BVE-BTP}, we see that for a small number of horizontal cells ($N_x=100$), the computational time required for both models is similar. However, using larger numbers of horizontal cells ($N_x=1000$), solving the BVE-model becomes distinctly faster than solving the BTP-model. Again, the extra time required by the BTP-model is used to solve a linear system for pressure, whose dimension grows by increasing the number of horizontal cells.

Example 2 shows a reduction in computational complexity of the BVE-model over that of the BTP-model. This is a consequence of the reduced number of unknowns.

\subsection{The BVE-Model and Saturation Overshoots}
\begin{figure}
\hspace{-0.4cm}
 \subfigure[VE-model]{
 \includegraphics[scale=0.21]{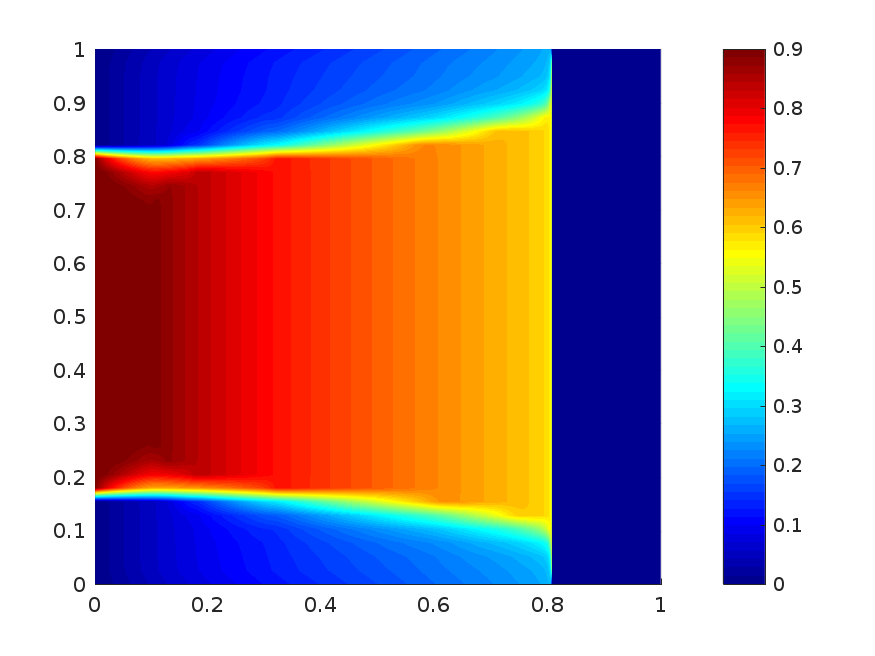}
 \label{subfig:VEmodel}
 }\hspace{-0.6cm}
  \subfigure[BVE-model]{
 \includegraphics[scale=0.21]{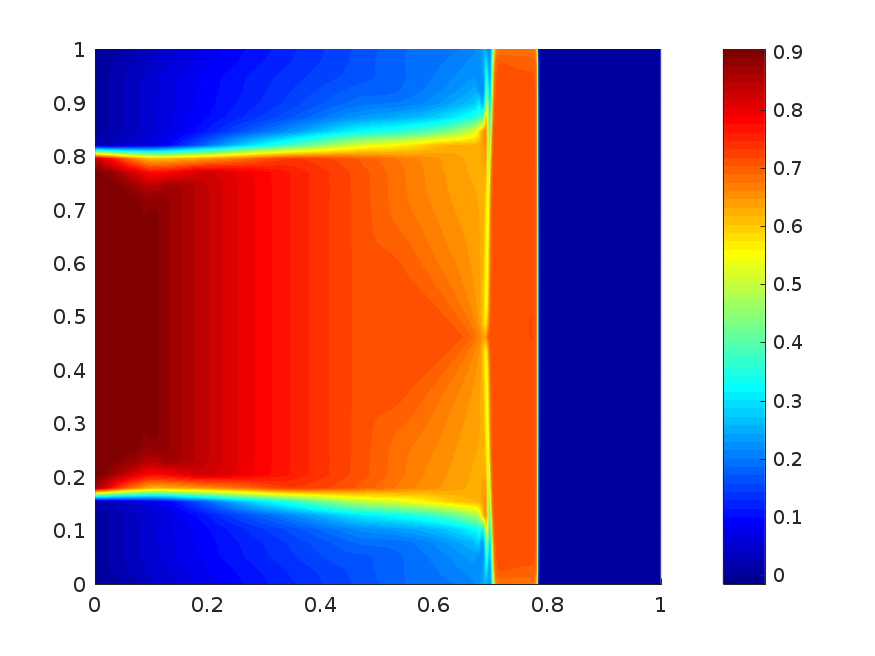}
  \label{subfig:BVEmodel}
 }\hspace{-0.6cm}
  \subfigure[Middle Layer]{
 \includegraphics[scale=0.15]{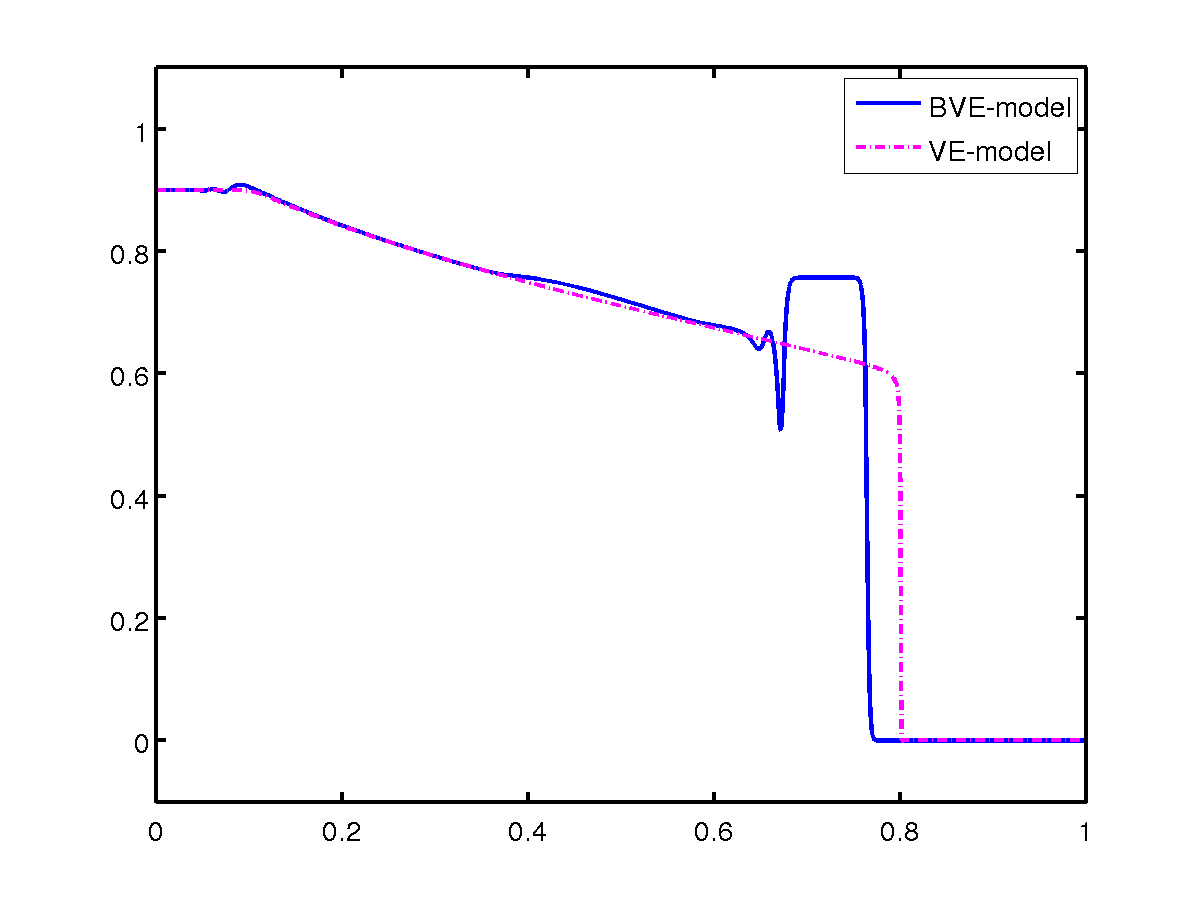}
  \label{subfig:1D-BVE-VE}
 }
 \caption{Numerical solution of the VE-model in (a), the BVE-model in (b) and saturation profile in the middle layer for both models in (c), using $M=2,\, \beta_x=\beta_z=10^{-6},\, N_x=2000$, $N_z=40$, and $T=0.6$. }
 \label{fig:BtwoLdifferentT}
\end{figure}

In this section, we emphasize the ability of the BVE-model \eqref{eq:Brinkman1}, \eqref{eq:Brinkman2} to drive saturation overshoots impacting the spreading speed. We set $\kappa=1$ and apply the inflow condition 
 \begin{align*}
  S_{\text{inflow}}(z)=\left\{ 
\begin{array}{c l l}
	0 \quad &: & z\leq \frac{1}{4} \text{ and } z>\frac{3}{4},\\
	0.9 \quad &: &\frac{1}{4}<z\leq \frac{3}{4}.
\end{array} \right.
 \end{align*}
In Figure \ref{subfig:VEmodel} we show the numerical solution for the VE-model and Figure \ref{subfig:BVEmodel} presents that of the BVE-model using a Cartesian grid with $N_x\times N_z=2000\times 40$, the parameters $M=2$, $T=0.6$, and $\beta_x=\beta_z=10^{-6}$.
 
Compared to the computations for the VE-model, Figure \ref{fig:BtwoLdifferentT} shows that the BVE-model describes saturation overshoots. Moreover, the spreading speed of the invading fluid $U\approx 1.27$ is smaller than that using the VE-model $u\approx 1.33$. Mathematically, the phenomenon of saturation overshoots is identified by undercompressive waves, which are known to be slower than classical compressive waves \cite{Kissling2015}. Yortsos and Salin \cite{YortsosSalin2006} expect that sharp bounds on the spreading speed might be found using these waves. 

Figure \ref{subfig:1D-BVE-VE} shows the saturation profile of both models in the middle layer of the domain, where saturation overshoots and its effect on the spreading speed using the BVE-model are more visible.

We conjecture that for $\beta_x,\,\beta_z,\epsilon_x,\,\epsilon_z\rightarrow 0$ the saturation overshoots in the BVE-model persist.

\section{Summary}
\label{sec:conclusion}

We studied two-phase flow for Darcy and Brinkman regimes. In Darcy regimes, we presented different models that utilize the extreme geometry of asymptotically flat domains to reduce the computational complexity of the two-phase flow model: the VI-model, the VE-model and the Multiscale model. We performed several numerical comparisons of the models and found out that the VE-model shows excellent agreement with the full two-phase flow model, while still being computationally more efficient. On the contrary, the VI-model shows limited accuracy. In addition to the numerical comparisons, we proved the equivalence between the VE-model and the Multiscale model.  

In flat domains of Brinkman type, we applied asymptotic analysis to the dimensionless BTP-model. We derived the BVE-model, which is a nonlocal evolution equation for saturation only. It involves a third-order term of mixed-time and spatial- derivatives. We showed by numerical examples that the BVE-model is a proper reduction of the full BTP-model in flat domains. The BVE-model has a similar structure as the model in \cite{HG1993}, for which dynamical capillarity effects drive overshooting. Here overshooting occurs within a completely different physical regime.


\begin{acknowledgements}
The first author would like to thank the DAAD for the financial support.  
\end{acknowledgements}

\bibliographystyle{plain}
\bibliography{ms}


%
%

\end{document}